\documentclass[12pt,reqno]{amsart}
\usepackage{latexsym}
\usepackage{amssymb}
\usepackage{mathrsfs}
\usepackage{amsmath}
\usepackage{color}
\usepackage{enumerate}
\usepackage[latin1]{inputenc}
\usepackage[colorlinks=true, linkcolor=blue, citecolor=blue]{hyperref}

\usepackage{color}

\usepackage[left=2.2cm,top=2.3cm,right=2.2cm]{geometry}

\def\1{\raisebox{2pt}{\rm{$\chi$}}}

\newtheorem{theorem}{Theorem}[section]

\newtheorem{lemma}[theorem]{Lemma}

\theoremstyle{definition}
\newtheorem{definition}[theorem]{Definition}
\newtheorem{remark}[theorem]{Remark}

\newcommand{\R}{{\mathbb R}}

\newcommand{\N}{{\mathbb N}}
\newcommand{\Z}{{\mathbb Z}}

\newcommand{\M}{{\mathcal M}}

\newcommand{\defeq}{\mathrel{\mathop:}=} 
\newcommand{\co}{\mskip0.5mu\colon\thinspace}   

\newcommand{\eps}{{\varepsilon}}
\def\1{\raisebox{2pt}{\rm{$\chi$}}}

\newcommand{\Lip}{\operatorname{Lip}}

\newcommand{\len}{\operatorname{len}} 
\def\vint_#1{\mathchoice%
        {\mathop{\kern 0.2em\vrule width 0.6em height 0.69678ex depth -0.58065ex
                \kern -0.8em \intop}\nolimits_{\kern -0.4em#1}}%
        {\mathop{\kern 0.1em\vrule width 0.5em height 0.69678ex depth -0.60387ex
                \kern -0.6em \intop}\nolimits_{#1}}%
        {\mathop{\kern 0.1em\vrule width 0.5em height 0.69678ex
            depth -0.60387ex
                \kern -0.6em \intop}\nolimits_{#1}}%
        {\mathop{\kern 0.1em\vrule width 0.5em height 0.69678ex depth -0.60387ex
                \kern -0.6em \intop}\nolimits_{#1}}}
\def\vintslides_#1{\mathchoice%
        {\mathop{\kern 0.1em\vrule width 0.5em height 0.697ex depth -0.581ex
                \kern -0.6em \intop}\nolimits_{\kern -0.4em#1}}%
        {\mathop{\kern 0.1em\vrule width 0.3em height 0.697ex depth -0.604ex
                \kern -0.4em \intop}\nolimits_{#1}}%
        {\mathop{\kern 0.1em\vrule width 0.3em height 0.697ex depth -0.604ex
                \kern -0.4em \intop}\nolimits_{#1}}%
        {\mathop{\kern 0.1em\vrule width 0.3em height 0.697ex depth -0.604ex
                \kern -0.4em \intop}\nolimits_{#1}}}

\newcommand{\aveint}[2]{\mathchoice%
        {\mathop{\kern 0.2em\vrule width 0.6em height 0.69678ex depth -0.58065ex
                \kern -0.8em \intop}\nolimits_{\kern -0.45em#1}^{#2}}%
        {\mathop{\kern 0.1em\vrule width 0.5em height 0.69678ex depth -0.60387ex
                \kern -0.6em \intop}\nolimits_{#1}^{#2}}%
        {\mathop{\kern 0.1em\vrule width 0.5em height 0.69678ex depth -0.60387ex
                \kern -0.6em \intop}\nolimits_{#1}^{#2}}%
        {\mathop{\kern 0.1em\vrule width 0.5em height 0.69678ex depth -0.60387ex
                \kern -0.6em \intop}\nolimits_{#1}^{#2}}}

\title[Self-improvement of weighted pointwise inequalities]{Self-improvement of weighted pointwise \\ inequalities on open sets}

\author[S.\! Eriksson-Bique]{Sylvester Eriksson-Bique}   
\address[S.E.-B.]{  Department of mathematics, UCLA, 520 Portola Plaza, Los Angeles CA 90095, USA }
\email{syerikss@math.ucla.edu}

\author[J.\! Lehrb\"ack]{Juha Lehrb\"ack}

\address[J.L.]{University of Jyvaskyla, Department of Mathematics and Statistics, P.O. Box 35, FI-40014 University of Jyvaskyla, Finland} \email{juha.lehrback@jyu.fi}

\author[A.V.\! V\"ah\"akangas]{Antti V. V\"ah\"akangas}
\address[A.V.V.]{University of Jyvaskyla, Department of Mathematics and Statistics, P.O. Box 35, FI-40014 University of Jyvaskyla, Finland} 
\email{antti.vahakangas@iki.fi}

\date{\today}

\begin{document}



\begin{abstract}
We prove a general self-improvement property for 
a family of weighted pointwise inequalities
on open sets, including 
pointwise Hardy inequalities with distance weights. 
For this purpose we introduce and study the classes of
$p$-Poincar\'e and $p$-Hardy weights for an open set $\Omega\subset X$,
where $X$ is a metric measure space. 
We also apply the self-improvement of weighted
pointwise Hardy inequalities in connection with
usual integral versions of Hardy inequalities.
\end{abstract}

\maketitle

\section{Introduction}

This paper is continuation of a general program related to various self-improving
phenomena, including Poincar\'e and Hardy inequalities and
uniform fatness; see e.g.~\cite{MR1869615,MR2415381,MR1948106,MR946438} for earlier results 
and \cite{SEB,MR3976590,MR3673660} for recent work by the authors.
In this paper we introduce a class of $p$-Hardy weights and consider for such weights $w$ 
the pointwise $(p,w)$-Hardy inequality
\begin{equation}\label{e.pw_intro}
\lvert u(x)\rvert \leq  C\, d(x,\Omega^c) 
\sup_{0<r<\kappa d(x,\Omega^c)} \biggl( \frac{1}{w(B(x,r))}\int_{B(x,r)} g(y)^p w(y)\,d\mu(y) \biggr)^{\frac{1}{p}}.
\end{equation}
Here $\Omega$ is an open subset of a metric space $X$,
$d(x,\Omega^c)$ denotes the distance from $x\in\Omega$ to the complement $\Omega^c=X\setminus\Omega$, and
$g$ is a (bounded) upper gradient of $u \in \Lip_{0}(\Omega)$; see Sections~\ref{s.notation}
and~\ref{s.weights} for definitions.
Our main result, Theorem~\ref{t.main_improvement}, 
shows that these inequalities are self-improving with respect to the
exponent $p$: if 
a pointwise $(p,w)$-Hardy inequality holds in $\Omega$ with an exponent $1<p<\infty$,
then, under suitable assumptions,  
there exists $1<q<p$ such that also a pointwise $(q,w)$-Hardy inequality holds in $\Omega$. 
The unweighted case $w=1$ corresponds to the pointwise $p$-Hardy inequality, for which
the self-improvement was proved in~\cite{MR3976590}. 
Our approach 
relies on the basic ideas and techniques developed in~\cite{SEB,MR3976590}. However,
unlike the self-improvement of pointwise $p$-Hardy inequalities, which was known
already before the work in~\cite{MR3976590} indirectly via the self-improvement of 
uniform $p$-fatness (see~\cite{MR1869615,MR946438}) and the equivalence between these two concepts 
(see~\cite{MR2854110}), the present self-improvement for the weighted pointwise $p$-Hardy inequalities is
previously unknown. In particular,
our main result is new even for $X=\R^n$, equipped with the Euclidean distance and the Lebesgue measure.

The self-improvement of the pointwise $(p,w)$-Hardy inequality and
a weighted maximal function theorem  
show that inequality~\eqref{e.pw_intro}, for every $x\in\Omega$, implies the 
integral version of the $(p,w)$-Hardy inequality, that is,
\begin{equation}\label{e.intro_w_hardy}
\int_\Omega\frac{\lvert u(x)\rvert^p}{d(x,\Omega^c)^p}\,w(x)\,d\mu(x)
\le C\int_\Omega g(x)^p w(x)\,d\mu(x);
\end{equation}
see Section~\ref{s.applications} for details.
This implication is not immediate
from inequality~\eqref{e.pw_intro}, since the maximal operator
is not typically bounded on $L^1(X)$. In some sense the inbuilt self-improvement
of pointwise Hardy inequalities provides a mechanism to bypass
the lack of the $L^1$-boundedness for the maximal operator.

An important model case of~\eqref{e.intro_w_hardy} is the weighted  $(p,\beta)$-Hardy inequality  in $\R^n$,  
with  $w(x)=d(x,\Omega^c)^\beta$, for $\beta\in\R$;
 see~\cite{Kufner1985,Necas1962}.
Corresponding pointwise theory was developed in~\cite{KoskelaLehrback2009}, 
but in order to be able to apply the maximal function theorem, 
it was necessary to assume {\it a priori} the validity of a stronger variant of~\eqref{e.pw_intro} in terms of an exponent $1< q<p$.
With the self-improvement
results of the present work, the starting point in the  weighted  
pointwise  Hardy inequalities as in~\cite{KoskelaLehrback2009}
can now be taken  to be  the natural candidate involving only the exponent $p$,
at least for $\beta\ge 0$. 
More motivation and explanation related to  (weighted)  pointwise Hardy inequalities 
in Euclidean spaces will be given in Section~\ref{s.applications}. 

Often the theory of weighted inequalities is concerned with doubling weights.
In the present setting the natural assumption is a weaker \emph{semilocal} doubling
condition with respect the open set  $\Omega\subsetneq X$.  This class of
weights is introduced in Section~\ref{s.weights}, 
where we also prove some technical lemmas for such weights.
As a tool in pointwise $(p,w)$-Hardy inequalities we also use a related class of 
$p$-Poincar\'e weights for $\Omega$, see Section~\ref{s.poinc_weight}. 
In Section~\ref{s.char} we define the $p$-Hardy weights, which will be crucial
for the pointwise $(p,w)$-Hardy inequalities, and in Section~\ref{s.self_improvement}
we establish a self-improvement result for $p$-Hardy weights. This plays a key role
also in the self-improvement of pointwise $(p,w)$-Hardy inequalities, since in
Section~\ref{s.pw_hardy} we show that $w$ being a $p$-Hardy weight is equivalent
to the validity of the pointwise $(p,w)$-Hardy inequality. Finally, Section~\ref{s.applications}
contains the applications related to integral versions of weighted Hardy inequalities.

\section{Notation and auxiliary results}\label{s.notation}

We make the standing assumption  
that $X=(X,d,\mu)$, with $\# X\ge 2$, 
is a metric measure space equipped with a metric $d$ and a 
positive complete $D$-doubling Borel regular
measure $\mu$ such that
$0<\mu(B)<\infty$ and
\begin{equation}\label{e.doubling}
\mu(2B) \le D\, \mu(B)
\end{equation}
for some $D>1$ and for all balls $B=B(x,r)=\{y\in X\mid d(y,x)<r\}$. 
Here we use for $0<\lambda <\infty$ the notation $\lambda B=B(x,\lambda r)$.
It follows that the 
space $X$ is separable (see e.g.\ \cite[Proposition 1.6]{MR2867756}) and 
$\mu(\{x\})=0$ for every $x\in X$ by \cite[Corollary 3.9]{MR2867756}.

For us, a {\em curve} is a rectifiable and continuous
mapping $\gamma \co [a,b] \to X$. By $\Gamma(X)$ we denote the set of all curves in $X$. The
length of a curve $\gamma\in \Gamma(X)$ is written as $\len(\gamma)$.
A curve $\gamma\colon [a,b]\to X$ {\em connects}
$x\in X$ to $y\in X$ (or a point $x\in X$ to a set $E\subset X$), if $\gamma(a)=x$ and $\gamma(b)=y$ ($\gamma(b)\in E$, respectively).
We assume throughout that the space $X$ is $C_{\mathrm{QC}}$-quasiconvex for some $C_{\mathrm{QC}}\ge 1$, that is,
for every $x,y\in X$ there exists a curve $\gamma$ connecting $x$ to $y$
such that $\len(\gamma)\le C_{\mathrm{QC}}d(x,y)$.
 
Fix $x,y \in X$, $E \subset X$ and $\nu \geq 1$. The collection $\Gamma(X)^\nu_{x,y}$ is the set of all 
curves that connect $x$ to $y$ and whose lengths are at most $\nu d(x,y)$. 
The set of all curves that connect $x$ to $E$ and whose lengths are at most $\nu d(x,E)$ is denoted by $\Gamma(X)^\nu_{x,E}$. 

A Borel function $g\ge 0$ on $X$ is an 
{\em upper gradient} of function $u\colon X \to \R$, if for all curves $\gamma\colon [a,b] \to X$, we have
\begin{equation}\label{e.modulus}
\lvert u(\gamma(a))-u(\gamma(b))\rvert \le \int_\gamma g\,ds.
\end{equation}

The space of Lipschitz functions on $X$ is denoted by $\Lip(X)$. By definition $u\in \Lip(X)$ if there
exists a constant $\lambda>0$ such that
\[
\lvert u(x)-u(y)\rvert \le \lambda d(x,y),\qquad \text{\ for all }\ x,y\in X.
\]
When $\Omega \subset X$ is an open set, we denote by 
$\Lip_{0}(\Omega)$
the space of all Lipschitz
functions on $X$ that vanish on $\Omega^c= X \setminus \Omega$. 
The set of lower semicontinuous functions on $X$ is denoted by $LC(X)$. 
 
Recall that 
\[
u_E=\vint_{E} u\,d\mu =\frac{1}{\mu(E)}\int_E u(y)\,d\mu(y)
\]
is the integral average of a function $u\in L^1(E)$ over a measurable set $E\subset X$
with $0<\mu(E)<\infty$. 
If $1\le p<\infty$ and $u\colon X\to \R$ is a $\mu$-measurable function, then $u\in L^p_{\textup{loc}}(X)$ 
means that for each $x_0\in X$ there exists $r>0$ such that $u\in L^p(B(x_0,r))$, 
that is, $\int_{B(x_0,r)} \lvert u\rvert^p\,d\mu<\infty$.
The characteristic function of a set $E\subset X$
is denoted by $\mathbf{1}_{E}$; that is, $\mathbf{1}_{E}(x)=1$ if $x\in E$
and $\mathbf{1}_{E}(x)=0$ if $x\in X\setminus E$.

\section{Weights and restricted maximal functions for open sets}\label{s.weights}

We need several classes of weights for open sets. 
To avoid pathological situations, we assume throughout the paper that
the open sets $\Omega\subset X$ under consideration are nonempty.

\begin{definition}\label{def.semilocal}
Let $\Omega\subset X$ 
be an open set.
A non-negative Borel function $w$ in $X$
is a weight for $\Omega$, if 
$\int_B w(x)\,dx<\infty$ for all balls $B\subset X$ and 
$w(x)>0$ for almost
every $x\in \Omega$. If $E\subset X$ is a measurable set, then we write $w(E)=\int_E w\,d\mu$.
\end{definition}

We impose the following localized doubling condition on the
weight $w$. We remark that there are also other
uses for the term \emph{semilocally doubling} in the literature,
see e.g.~\cite{BBsemilocal}. 
In our definition ``local'' refers to the fact that the condition is
required only for points $x\in\Omega$, but ``semi'' is added since the
balls need not be contained in $\Omega$.

\begin{definition}\label{def.semilocal_dbl}
Let $\Omega\subsetneq X$ be an open set
and let $w$ be a weight for $\Omega$. We say that $w$ 
is semilocally doubling for $\Omega$
if for every $\kappa>0$ there exists a constant $D(w,\kappa)\ge 1$ such that
\[
0<w(B(x,r))\le D(w,\kappa)w(B(x,r/2))<\infty
\]
for all $x\in \Omega$ and $0<r\le \kappa d(x,\Omega^c)$.
\end{definition}

In some of our results we will need the following regularity property of $w$.

\begin{lemma}\label{lem:outreg} 
Let  $\Omega\subset X$ be an open set
and let $w$ be a weight for $\Omega$. Then $w$ is outer regular,
that is,
for every Borel set $E\subset X$ and every $\varepsilon>0$, there
exists an open set $V\supset E$ such that
$w(V)\le w(E)+\varepsilon$.
\end{lemma}

\begin{proof}
Let $\overline{X}$ be a completion of $X$. We remark that $X$ could fail to be a Borel subset of its completion. 
We denote by $\mathscr{B}(X)$ and $\mathscr{B}(\overline{X})$ the Borel
sets of $X$ and $\overline{X}$, respectively.
The measures $\mu$ and $w\mu$ extend to Borel regular measures $\overline{\mu}$ and $\overline{w\mu}$ on $\overline{X}$, 
and $\overline{\mu}$ is doubling, by \cite[Lemma 1]{MR1678044}.
More precisely
\begin{equation}\label{e.saksman}
\{F\in\mathscr{B}(\overline{X})\mid F\cap X\in\mathscr{B}(X)\}=\mathscr{B}(\overline{X}),
\end{equation} 
and therefore one can define 
$\overline{\mu}(F) = \mu(F \cap X)$ and $\overline{w\mu}(F) = w(F \cap X)$ 
for each $F\in\mathscr{B}(\overline{X})$; see the proof of \cite[Lemma 1]{MR1678044}.
This defines the extended measures as Borel measures that are finite on balls,
and the Borel regular (complete) extended measures are obtained by completion. 
The space $\overline{X}$ is complete and the measure $\overline{\mu}$ doubling; thus $\overline{X}$ is proper
by \cite[Proposition 3.1]{MR2867756}.  
Hence, the measure $\overline{w\mu}$ is outer regular on $\overline{X}$ by \cite[Theorem 7.8]{MR767633}.

Let $E\in\mathscr{B}(X)$ and $\eps>0$. By using $\sigma$-algebra arguments, one can
show that $E=F\cap X$ for some $F\in\mathscr{B}(\overline{X})$.
By the outer regularity of $\overline{w\mu}$, there
exists an open set $U$ in $\overline{X}$ such that
$U\supset F$ and $\overline{w\mu}(U)\le \overline{w\mu}(F)+\eps$.
We define $V=U\cap X$, which is an open subset of $X$. Then $V\supset E$ and
\[
w(V)=w(U\cap X)=\overline{w\mu}(U)\le \overline{w\mu}(F)+\eps=w(F\cap X)+\eps=w(E)+\eps.
\]
This shows that $w$ is outer regular.
\end{proof}

Let $\Omega\subsetneq X$ be an open set and fix
a weight $w$ for $\Omega$.
Let $0< \kappa  <\infty$ and $1\le p<\infty$, and let $f$ be a measurable
function in $X$. We define  restricted weighted maximal functions $\M_{p,w,\kappa} f$ and $\M_{p,w,\kappa}^R f$ at $x\in \Omega$ by
\[
 \M_{p,w,\kappa} f(x) \defeq 
\sup_{0<r<\kappa d(x,\Omega^c)} \left( \frac{1}{w(B(x,r))}\int_{B(x,r)} \lvert f\rvert^p w\,d\mu\right)^\frac{1}{p}
\]
and
\[
\M_{p,w,\kappa}^R f(x)
\defeq \sup_{0<r<\min \{\kappa d(x,\Omega^c),R\}} \left( \frac{1}{w(B(x,r))}\int_{B(x,r)} \lvert f\rvert^p w\,d\mu\right)^\frac{1}{p}.
\]
Observe that $0<w(B(x,r))<\infty$ for all balls $B(x,r)$
that appear within the supremums.
The maximal functions  $\M_{p,w,\kappa} f$ and $\M_{p,w,\kappa}^R f$
are lower semicontinuous in $\Omega$. This
follows easily using monotone convergence theorem and the fact that 
$B=\bigcup_{0<\eps<1} (1-\eps)B$
for all balls $B\subset X$.

The following lemmas are adaptations of similar results from
our prior work \cite{SEB, MR3976590}. 
Although the methods are the same, 
we provide here the full proofs
due to subtle technical differences.

\begin{lemma}\label{lem:HLmaxmaxest}
Suppose that $w$ is a semilocally doubling weight for an open set  $\Omega\subsetneq X$. 
Assume that $1\le q<\infty$ and $\kappa> 1$, and let 
$f \in L^q_{\mathrm{loc}}(X)$,
$x\in \Omega$
and $\tau>0$ be such that
\[
\M_{q,w,2\kappa} f(x)\le \tau.
\]
Fix $\Lambda>0$ and define  
\[
E_{\Lambda} =\{y\in \Omega \mid \M_{q,w,2\kappa}^{\kappa d(x,\Omega^c)} f(y) > \Lambda \tau\}.
\]
Then 
\begin{equation}\label{eq:localMaxMax}
\M_{1,w,\kappa} \mathbf{1}_{E_{\Lambda}} (x) \leq \frac{D(w,10 \kappa)^4}{\Lambda^q}.
\end{equation} 
\end{lemma}

\begin{proof}
Fix $0<r<\kappa d(x,\Omega^c)$ and let
$B=B(x,r)$.
We need to show that
\begin{equation}\label{e.M_single_ball}
\frac{1}{w(B)}\int_{B} \mathbf{1}_{E_{\Lambda}}w\, d\mu
\le \frac{D(w,10\kappa)^4}{\Lambda^q}.
\end{equation}
The proof of \eqref{e.M_single_ball} uses a covering argument.
For each $y\in E_{\Lambda}\cap B$ we fix
a ball $B_y=B(y,r_y)$ of radius $0<r_y<\min\{2\kappa d(y,\Omega^c),\kappa d(x,\Omega^c)\}$ such that 
\begin{equation}\label{e.sel}
\biggl(\frac{1}{w(B_y)}\int_{B_y}   \lvert f\rvert^q w\,d\mu\biggr)^{\frac{1}{q}} >\Lambda \tau.
\end{equation}
There are two cases to consider.

\vskip1mm
\textbf{Case 1: There exists $y\in E_{\Lambda}\cap B$ with $r<r_y$.} Then
$B(y,r_y)\subset B(x,2r_y)\subset B(y,3r_y)$ and $3r_y\le 6\kappa d(y,\Omega^c)$. By semilocal doubling, we have
$w(B(y,3r_y))\le  D(w,10\kappa)^2 w(B(y,r_y))$. Observe that $2r_y<2\kappa d(x,\Omega^c)$.
Therefore,
\begin{align*}
\frac{1}{w(B)}\int_{B} \mathbf{1}_{E_{\Lambda}} w\,d\mu
&\le 1< \frac{\frac{1}{w(B_y)}\int_{B_y}   \lvert f\rvert^q w\,d\mu}{\Lambda^q \tau^q}
\\&\le  \frac{D(w,10\kappa)^2\frac{1}{w(B(x,2r_y))}\int_{B(x,2r_y)}   \lvert f\rvert^qw \,d\mu}{\Lambda^q \tau^q}
\\&\le \frac{D(w,10\kappa)^2(\mathcal{M}_{q,w,2\kappa} f(x))^q}{\Lambda^q\tau^q}
\le \frac{D(w,10\kappa)^4}{\Lambda^q},
\end{align*}
proving inequality \eqref{e.M_single_ball}. 

\vskip1mm
\textbf{Case 2: For each $y\in E_{\Lambda}\cap B$ we have $r\ge r_y$.} The $5r$-covering lemma 
\cite[Lemma 1.7]{MR2867756} yields a pairwise disjoint subcollection $\mathcal{B}\subset \{B_y \mid y\in E_{\Lambda}\cap B\}$ 
of balls such that $E_{\Lambda}\cap B\subset \bigcup_{B'\in\mathcal{B}} 5B'$.
Hence, by \eqref{e.sel} and the fact that $5r_y\le 10 \kappa d(y,\Omega^c)$ for
every $y\in E_\Lambda \cap B$,
\begin{align*}
\frac{1}{w(B)}\int_{B} \mathbf{1}_{E_{\Lambda}}w\, d\mu
&\le \frac{1}{w(B)}\sum_{B'\in\mathcal{B}} w(5B')\\
&\le \frac{D(w,10\kappa)^3}{w(B)}\sum_{B'\in\mathcal{B}} w(B')\\
&\le \frac{D(w,10\kappa)^3}{\Lambda^q\tau^q w(B)}\sum_{B'\in\mathcal{B}}
\int_{B'}  \lvert f\rvert^qw \,d\mu.
\end{align*}
Since $r_{B'}\le r$, we have 
$B'\subset 2B=B(x, 2r)$ for every $B'\in\mathcal{B}$. Also,
since $2r< 2\kappa d(x,\Omega^c)$, we have $w(2B) \leq D(w,10\kappa) w(B)$.
Consequently, 
inequality \eqref{e.M_single_ball} follows from the 
estimates
\begin{align*}
\frac{1}{w(B)}\int_{B} \mathbf{1}_{E_{\Lambda}}w\, d\mu
&\le \frac{D(w,10\kappa)^4}{\Lambda^q\tau^qw(2B)}\int_{2B}   \lvert f\rvert^qw \,d\mu\\
&\le \frac{D(w,10\kappa)^4(\mathcal{M}_{q,w,2\kappa}   f (x))^q}{\Lambda^q\tau^q}\le  \frac{D(w,10\kappa)^4}{\Lambda^q}. \qedhere
\end{align*}
\end{proof}

The next 
approximation lemma is a variant of \cite[Lemma 3.7]{SEB}. The outer regularity of the weight, see Lemma \ref{lem:outreg}, is needed in the proof.
Recall that a Borel function $g\colon X\to [0,\infty)$ is {\em simple}, if it can be expressed as
$g=\sum_{j=1}^k a_j\mathbf{1}_{E_j}$ for
some real numbers  $a_j>0$ and Borel sets $E_j\subset X$, $j=1,\ldots,k$.

\begin{lemma}\label{lem:approx}   
Suppose that $w$ is a semilocally doubling weight for an open set $\Omega$.
Assume that $1\le p<\infty$ and $\kappa> 1$, and let 
$g\colon X\to [0,\infty)$ be a simple Borel function. Then, for each finite set $F\subset \Omega$ and every $\varepsilon>0$, 
there exists a non-negative and bounded $g_{F,\varepsilon}\in LC(X)$ such that
$g(y)\le g_{F,\varepsilon}(y)$ for all  $y\in X\setminus F$  and
$\M_{p,w,\kappa}g_{F,\varepsilon}(x) \leq \M_{p,w,\kappa}g(x)+\varepsilon$ for every $x\in F$.
\end{lemma}

\begin{proof} 
It suffices to prove the claim for singletons $F=\{x\}$, 
since for $F = \{x_1, \dots, x_n\}$ the function $g_{F,\varepsilon}$ can be obtained as the minimum 
of the functions $g_{\{x_i\}, \varepsilon}$. Fix $x\in \Omega$ and $\varepsilon>0$.

\vskip1mm
\textbf{Step 1: proving the claim for $g=\mathbf{1}_E$ with a Borel set $E$.} 
We show that there exists an open set $U\subset X$ such that $\mathbf{1}_E\le \mathbf{1}_U$ in $X\setminus \{x\}$ and
\begin{equation}\label{e.verysimple}
\M_{p,w,\kappa}(\mathbf{1}_U-\mathbf{1}_E)(x)<\varepsilon.
\end{equation}
For each $m\in \Z$, we set 
\[A_m=\{y\in X \mid 2^{m-1}<d(x,y)<2^{m+1}\}.\]
Observe that each $y\in X$ belongs to at most two annuli $A_m$. 
Moreover, if $m\in\Z$ then by outer regularity of the weight $w$ (Lemma \ref{lem:outreg})
and the fact that $A_m$ is open, there is an open set $U_m\subset A_m$
such that
\begin{equation}\label{e.measures}
\begin{split}
&A_m\cap E\subset U_m\quad \text{\ and }\ \quad \\ &w(U_m\setminus E)=w(U_m\setminus (A_m\cap E))\le   \frac{\eps^pw(A_m)}{2D(w,4\kappa)^2} .
\end{split}
\end{equation}
In the case $w(A_m)=0$ we can choose $U_m=A_m$. 
Define $U=\bigcup_{m\in\Z} U_m$. Then
\begin{equation}\label{e.inclusion}
E\setminus \{x\}\subset \bigcup_{m\in\Z} (A_m\cap E)\subset \bigcup_{m\in\Z} U_m=U.
\end{equation}
As a consequence, we have $\mathbf{1}_E(y)\le \mathbf{1}_U(y)$ for every $y\in X\setminus \{x\}$. 

To prove \eqref{e.verysimple}, we let $B(x,r)\subset X$ be a ball 
with $0<r<\kappa d(x,\Omega^c)$.  Then
$\mathbf{1}_U-\mathbf{1}_E=\mathbf{1}_{U\setminus E}$ $\mu$-almost everywhere, and therefore
by \eqref{e.measures} we obtain 
\begin{align*}
\frac{1}{w(B(x,r))}\int_{B(x,r)} \lvert \mathbf{1}_U-\mathbf{1}_E\rvert^pw\,d\mu&=\frac{1}{w(B(x,r))}\int_{B(x,r)} \mathbf{1}_{U\setminus E}w\,d\mu\\
&\le \frac{1}{w(B(x,r))} \int_X\sum_{m=-\infty}^{\lceil \log_2 r\rceil} \mathbf{1}_{U_m\setminus E}w\,d\mu\\
&= \frac{\varepsilon^{p}}{2D(w,4\kappa)^2w(B(x,r))}\sum_{m=-\infty}^{\lceil \log_2 r\rceil} w(A_m)\\
&\le \frac{\varepsilon^{p}}{D(w,4\kappa)^2}\frac{w(B(x,4r))}{w(B(x,r))}\le \eps^p\frac{w(B(x,r))}{w(B(x,r))}= \varepsilon^{p}.
\end{align*}
Inequality \eqref{e.verysimple} follows by 
raising this estimate to power $\frac{1}{p}$ and then taking
supremum over all balls $B(x,r)$ as above.

\vskip1mm
\textbf{Step 2: proving the claim for a simple Borel function $g=\sum_{j=1}^k a_j\mathbf{1}_{E_j}$.} 
By Step~1,
for each $j=1,\ldots,k$, there exists a non-negative and bounded 
$g_{\{x\},\varepsilon,j}\in  LC(X)$ such that $\mathbf{1}_{E_j}\le g_{\{x\},\varepsilon,j}$ in $X\setminus \{x\}$
and 
\begin{equation}\label{e.error}
\M_{p,w,\kappa}(g_{\{x\},\varepsilon,j}-\mathbf{1}_{E_j})(x) \le\frac{\varepsilon}{k\max_j a_j }.
\end{equation}
Define  $g_{\{x\},\varepsilon}=\sum_{j=1}^k a_jg_{\{x\},\varepsilon,j}$.
Then $g\le g_{\{x\},\varepsilon}$ in $X\setminus \{x\}$, 
and  by using the subadditivity and positive homogeneity of the maximal function
and inequality \eqref{e.error}, we  conclude that 
\begin{align*}
\M_{p,w,\kappa}g_{\{x\},\varepsilon}(x)&=\M_{p,w,\kappa}(g+g_{\{x\},\varepsilon}-g)(x)\\&\le \M_{p,w,\kappa}g(x)
+\M_{p,w,\kappa}(g_{\{x\},\varepsilon}-g)(x)\\&\le \M_{p,w,\kappa}g(x)
 + \sum_{j=1}^k a_j\M_{p,w,\kappa}(g_{\{x\},\varepsilon,j}-\mathbf{1}_{E_j})(x)\\&\le \M_{p,w,\kappa}g(x)+\varepsilon. \qedhere
\end{align*}
\end{proof}

\section{Local Poincar\'e inequalities in open sets}\label{s.poinc_weight}

In the sequel, we will need to assume that a suitable pointwise  Poincar\'e inequality holds
with respect to the weight $w$.

\begin{definition}\label{def:PIchar}
Let $1\le p<\infty$, let  $\Omega\subsetneq X$  be an open set and let $w$ be a weight for $\Omega$.
We say that $w$ is a $p$-Poincar\'e weight for $\Omega$, if
there are constants $C_{\mathrm{A}}>0$, $\nu>C_{\mathrm{QC}}$ and $\kappa> 1$ such that  
for each non-negative and bounded  $g \in LC(X)$ and every  $x,y \in \Omega$ with  
\[
d(x,y)< d(x,\Omega^c)/(3\kappa),
\] 
it holds that
\begin{equation}\label{eq:PIchar}
\inf_{\gamma \in \Gamma(X)^\nu_{x,y}} \int_{\gamma} g \,ds 
\leq   C_{\mathrm{A}}  \,d(x,y) \left( \M_{p,w,\kappa}^{\kappa d(x,y)} g (x) + \M_{p,w,\kappa}^{\kappa d(x,y)} g (y)\right).
\end{equation}
\end{definition}

Definition \ref{def:PIchar} for a $p$-Poincar\'e weight is slightly technical, since
it is adjusted to our later purposes. 
The following lemma provides a more familiar
variant of a $p$-Poincar\'e inequality that is sufficient for \eqref{eq:PIchar}. We emphasize
the local nature of these Poincar\'e inequalities with respect to $\Omega$; for instance, we only require inequality \eqref{e.familiar}
for balls $B$ satisfying $2\lambda B\subset \Omega$.
Compare also to~\cite{Hajlasz2002}, and the references therein, 
concerning Poincar\'e inequalities and
pointwise inequalities related to~\eqref{eq:PIchar}.

We write  
$u_{B;w}=\frac{1}{w(B)}\int_B u(x)w(x)\,d\mu(x)$ whenever $uw\in L^1(B)$ and $B$ is a ball in $X$.

\begin{lemma}\label{l.poinc_suff}
Let $1\le p<\infty$ and $1\le\lambda <\infty$,
let $\Omega\subsetneq X$  be an open set, and let
$w$ be a semilocally doubling weight for $\Omega$.
Suppose there
exists a constant $C_1$ such that for each $u\in \Lip(X)$ and for every bounded upper gradient $g$ of $u$ we have
\begin{equation}\label{e.familiar}
\frac{1}{w(B)}\int_{B}\lvert u-u_{B;w}\rvert w\,d\mu
 \le C_1r\left(\frac{1}{w(\lambda B)}\int_{\lambda B} g^p w\,d\mu\right)^{\frac{1}{p}},
\end{equation}
whenever $B=B(x,r)$ is a ball with
$2\lambda B\subset\Omega$. Then $w$ is a $p$-Poincar\'e weight for $\Omega$.
\end{lemma}


\begin{proof}
The proof has two steps.
\vskip1mm
\textbf{Step 1: We show that  
there exist constants $C_2=6C_1D(w,2^{-1})^2$  and $\kappa=3\lambda$ such that
\begin{equation}\label{e.valipala}
\lvert u(x)-u(y)\rvert
\le C_2d(x,y)
\bigl(\M_{p,w,\kappa}^{\kappa d(x,y)} g (x) + \M_{p,w,\kappa}^{\kappa d(x,y)} g (y)\bigr)
\end{equation}
for every $x,y\in \Omega$ with
$d(x,y)<d(x,\Omega^c)/(3\kappa)$.} Here $u$ and $g$ are as in the assumptions of the lemma.

Fix $x,y\in \Omega$, with $x\not=y$ and $r=d(x,y)<d(x,\Omega^c)/(9\lambda)$. Write 
$B_i=B(x,2^{-i}r)$, for every $i\in \N_0$.
A telescoping argument yields
\begin{align*}
\lvert u(x)-u_{B(x,r);w}\rvert
&\le\sum_{i=0}^\infty\lvert u_{B_{i+1};w}-u_{B_i;w}\rvert
\\
&\le\sum_{i=0}^\infty\frac{w(B_i)}{w( B_{i+1})}
\frac{1}{w(B_i)}\int_{B_i}\lvert u-u_{B_i;w}\rvert w\,d\mu
\\
&\le C_1 D(w,2^{-1})\sum_{i=0}^{\infty}(2^{-i}r)
\left(\frac{1}{w(\lambda B_i)}\int_{\lambda B_i} g^p w\,d\mu\right)^{\frac{1}{p}}
\\
&\le 2C_1 D(w,2^{-1})d(x,y) \M_{p,w,\kappa}^{\kappa d(x,y)} g (x).
\end{align*}
Observe that $B(x,r)\subset B(y,2r)$
and
$2r=2d(x,y)<d(x,\Omega^c)/(4\lambda)\le  d(y,\Omega^c)/(2\lambda)$. Thus, a similar telescoping argument gives 
\[
\lvert u(y)-u_{B(y,2r);w}\rvert\le 4C_1 D(w,2^{-1})d(x,y) \M_{p,w,\kappa}^{\kappa d(x,y)} g (y).
\]
Since $B(x,r)\subset B(y,2r)\subset B(x,4r)$, we also have
\begin{align*}
\lvert u_{B(x,r);w}-u_{B(y,2r);w}\rvert&\le \frac{1}{w(B(x,r))}\int_{B(x,r)}\lvert u-u_{B(y,2r);w}\rvert w\,d\mu
\\
&\le \frac{w(B(x,4r))}{w(B(x,r))} \frac{1}{w(B(y,2r))}\int_{B(y,2r)}\lvert u-u_{B(y,2r);w}\rvert w\,d\mu\\
&\le 2C_1D(w,2^{-1})^2 d(x,y) \left(\frac{1}{w(B(y,2\lambda r))}
\int_{B(y,2\lambda r)} g^p w\,d\mu\right)^{\frac{1}{p}}\\
&\le 2C_1 D(w,2^{-1})^2d(x,y) \M_{p,w,\kappa}^{\kappa d(x,y)} g (y).
\end{align*}
By combining the estimates above we obtain
\begin{align*}
\lvert u(x)-u(y)\rvert
&\le\lvert u(x)-u_{B(x,r);w}\rvert
+\lvert u_{B(x,r);w}-u_{B(y,2r);w}\rvert+\lvert u(y)-u_{B(y,2r);w}\rvert
\\
&\le 6C_1D(w,2^{-1})^2 d(x,y)
\bigl(\M_{p,w,\kappa}^{\kappa d(x,y)} g (x) + \M_{p,w,\kappa}^{\kappa d(x,y)} g (y)\bigr),
\end{align*}
and this completes the proof of inequality \eqref{e.valipala}.

\vskip1mm
\textbf{Step  2:  With the aid of inequality \eqref{e.valipala}, we show that $w$ is a $p$-Poincar\'e weight for $\Omega$.}
Let $g\in LC(X)$ be a non-negative and bounded function. 
Fix $x,y\in \Omega$ such that $0<d(x,y)<d(x,\Omega^c)/(3\kappa)$ and let $\delta>0$; 
here $\kappa=3\lambda$ by Step 1.
Define 
$u\colon X\to [0,\infty)$ by setting  
\begin{equation}\label{e.poinc_u_def}
  u(z)   = \inf_{\gamma} \int_\gamma h \,ds,\qquad   z\in X,  
\end{equation}
where 
\[
h = g + \bigl(\M_{p,w,\kappa}^{\kappa d(x,y)} g (x) + \M_{p,w,\kappa}^{\kappa d(x,y)} g (y)+\delta\bigr)
\] 
and the infimum is taken over all curves
$\gamma$ in $X$ connecting $z$ to $y$. 
Note that $h$ is a non-negative bounded Borel function, and  
clearly $u(y)=0$.
Fix $z_1,z_2\in X$ and consider any curve $\sigma$ connecting $z_1$ to $z_2$. 
We claim that
\begin{equation}\label{e.poinc_des}
\lvert u(z_1)-u(z_2)\rvert\le \int_\sigma h\,ds.
\end{equation}
From this it follows, in particular, that $h$ is an upper gradient of $u$. 
Moreover, since $X$ is quasiconvex and $h$ is bounded, estimate \eqref{e.poinc_des} 
implies that $u\in \Lip(X)$.

In order to prove  \eqref{e.poinc_des}, we may
assume that $u(z_1)>u(z_2)$. Fix $\eps>0$ and let $\gamma$ be a curve in $X$ that
connects $z_2$ to $y$ and satisfies inequality  
\[
u(z_2)\ge \int_\gamma h\,ds - \eps.
\]
Let $\sigma \gamma$ be the concatenation of $\sigma$ and $\gamma$. Then
\begin{align*}
\lvert u(z_1)-u(z_2)\rvert &=u(z_1)-u(z_2)\\
&\le \int_{ \sigma\gamma } h\,ds - \int_\gamma h\,ds + \eps=\int_{\sigma} h\,ds+\eps.
\end{align*}
The desired inequality \eqref{e.poinc_des}
follows by taking $\eps\to 0_+$.

Application of inequality \eqref{e.valipala}
to $u\in \Lip(X)$ and its bounded upper gradient $h$ gives
\[
\lvert u(x)-u(y)\rvert
\le C_2d(x,y)
\bigl(\M_{p,w,\kappa}^{\kappa d(x,y)} h (x) + \M_{p,w,\kappa}^{\kappa d(x,y)} h (y)\bigr)<\infty.
\]
Since $u(x)\ge \delta d(x,y)>0$ and $u(y)=0$, by \eqref{e.poinc_u_def} there is a curve $\gamma$ in $X$ 
connecting $x$ to $y$ such that  
\begin{equation}\label{e.poinc_strong}
\begin{split}
\int_\gamma g\,ds \ +\, &\left(\M_{p,w,\kappa}^{\kappa d(x,y)} g (x) + \M_{p,w,\kappa}^{\kappa d(x,y)} g (y)+\delta\right)\len(\gamma)
 \\&\qquad = \int_\gamma h\,ds \le 2u(x)
 = 2\lvert u(x)-u(y)\rvert \\&
 \qquad \le 2C_2d(x,y)
\bigl(\M_{p,w,\kappa}^{\kappa d(x,y)} h (x) + \M_{p,w,\kappa}^{\kappa d(x,y)} h (y)\bigr)\\&\qquad 
\le 2C_{2}\,d(x,y)\bigl(3\M_{p,w,\kappa}^{\kappa d(x,y)} g (x) + 3\M_{p,w,\kappa}^{\kappa d(x,y)} g (y)+2\delta\bigr)\\&\qquad
\le 6C_{2}\,d(x,y)\bigl(\M_{p,w,\kappa}^{\kappa d(x,y)} g (x) + \M_{p,w,\kappa}^{\kappa d(x,y)} g (y)+\delta\bigr).
\end{split}
\end{equation}
The second last inequality  follows from the sublinearity of maximal function and definition of $h$. 
From \eqref{e.poinc_strong} we see that $\len(\gamma)\le 6C_{2}\,d(x,y)$.
By taking $\delta\to 0_+$, we also obtain from \eqref{e.poinc_strong} that inequality
\eqref{eq:PIchar} holds, that is,
\[
\inf_{\gamma \in \Gamma(X)^\nu_{x,y}} \int_{\gamma} g \,ds 
\leq   C_{\mathrm{A}}  \,d(x,y) \left( \M_{p,w,\kappa}^{\kappa d(x,y)} g (x) + \M_{p,w,\kappa}^{\kappa d(x,y)} g (y)\right),
\]
with $C_{\mathrm{A}}=6C_2\,$, $\kappa=3\lambda$ and $\nu>\max\{C_{\mathrm{QC}},6C_2\}$.
\end{proof}

\section{The class of $p$-Hardy weights}\label{s.char}

The following class of weights turns out to be natural
in connection with pointwise Hardy inequalities; see Lemma~\ref{lem:HardyChar},
and compare also to the definition of $p$-Poincar\'e weights in Definition~\ref{def:PIchar}.

\begin{definition}\label{def:HardyChar} 
Let $1\le p<\infty$, let  $\Omega\subsetneq X$  be an open set, and let $w$ be a weight
for $\Omega$. We say that 
$w$ is a $p$-Hardy weight for $\Omega$ if there are constants $C_\Gamma>0$, $\nu >C_{\mathrm{QC}}$ and $\kappa> 1 $ 
such that for each non-negative and bounded $g \in LC(X) $ and every $x \in \Omega$, we have
\begin{equation}\label{eq:HardyChar} 
\inf_{\gamma \in \Gamma(X)^\nu_{x,\Omega^c}} \int_{\gamma} g \,ds \leq C_\Gamma \,d(x,\Omega^c)  \M_{p,w,\kappa} g (x).
\end{equation}
\end{definition}

Next we define  a convenient albeit slightly abstract $\alpha$-function that condenses the 
$p$-Hardy weight property, specifically inequality \eqref{eq:HardyChar}, in a single function. 
Indeed, despite the
complex appearance this function is
a very useful tool in the proof of the self-improvement for $p$-Hardy weight property.
 
\begin{definition}
Let $\Omega\subsetneq X$  be an open set and let $w$ be a weight for $\Omega$.
If $\tau\ge 0$, $\kappa>1$, $1\le p<\infty$ and $x\in \Omega$, we write
\[
\mathcal{E}^{\kappa,\tau}_{p,w,x,\Omega} = \{g \in   LC(X)   \mid  \M_{p,  w,\kappa} g (x) \leq \tau \text{ and }g(y) \in [0,1]\text{ for all }y\in X\}.
\]
If also $\nu>C_{\mathrm{QC}}$, then we write   
\begin{equation}\label{eq:alphaDef}
\alpha_{p,w,\Omega}(  \nu,\kappa,\tau) \defeq \sup_{x \in \Omega} \sup_{g \in \mathcal{E}^{\kappa,\tau}_{p,w,x,\Omega}} \frac{ \inf_{\gamma \in \Gamma(X)^\nu_{x,\Omega^c}}\int_\gamma g \,ds}{d(x,\Omega^c)}. 
\end{equation}
\end{definition}

The parameter $\nu$ is related to the maximum length of the curves $\gamma$,
since $\len(\gamma)\le \nu d(x,\Omega^c)$. The remaining parameters $\kappa$ and $\tau$ are used to control 
the non-locality and size, or ``level'', of the maximal function $\M_{p,w,\kappa} g(x)$.
%

The following lemma codifies the relationship between
inequality \eqref{eq:HardyChar} and the $\alpha$-function.

\begin{lemma} \label{lem:rewrite}
Let $\Omega\subsetneq X$ be an open set and let $w$ be a weight for $\Omega$. 
Assume that $\kappa>1$, $1\le p<\infty$ and $\nu>C_{\mathrm{QC}}$, and let 
$g\in LC(X)$ be such that $g(y)\in [0,1]$ for
every $y\in X$.
Then, for every $x\in \Omega$, we have 
 \begin{equation}
\inf_{\gamma \in \Gamma(X)^\nu_{x,\Omega^c}} \int_{\gamma} g \,ds 
\leq d(x,\Omega^c) \alpha_{p,w,\Omega}\bigl(\nu,\kappa,\left( \M_{p,w,\kappa} g (x) \right)\bigr). 
\end{equation}
\end{lemma}

\begin{proof}
Take any $g\in LC(X)$ with $g(y)\in [0,1]$ for all $y\in X$. Fix $x \in \Omega$
and write
\[\tau=\M_{p,w,\kappa} g (x) \ge 0 .\]
Then $g\in \mathcal{E}^{\kappa,\tau}_{p,w,x,\Omega}$, and by the definition of $\alpha_{p,w,\Omega}$ 
\[
\frac{\inf_{\gamma \in \Gamma(X)^\nu_{x,\Omega^c}} \int_{\gamma} g \,ds}{d(x,\Omega^c)}\le
\sup_{h \in \mathcal{E}^{\kappa,\tau}_{p,w,x,\Omega}} \frac{ \inf_{\gamma \in \Gamma(X)^\nu_{x,\Omega^c}}\int_\gamma h \,ds}{d(x, \Omega^c)}
\le \alpha_{p,w,\Omega}(  \nu,\kappa,\tau).
\]
The last step holds, since $x\in \Omega$.  
\end{proof}

In particular, from Lemma \ref{lem:rewrite} we
obtain the following sufficient condition for
$p$-Hardy weights in terms of a $\tau$-linear upped bound for the $\alpha$-function.

\begin{lemma}\label{lem:alphaClass}
Let $1\le p<\infty$, let $\Omega\subsetneq X$  be an open set and let $w$ be a weight for $\Omega$.
Suppose that there
are constants $\nu >C_{\mathrm{QC}}$, $\kappa> 1$ and $C_\alpha >0$ such that,
for any $\tau\ge 0$, we have
\[
\alpha_{p,w,\Omega}(\nu,\kappa,\tau)\le C_\alpha \tau.
\]
Then $w$ is a $p$-Hardy weight for $\Omega$.
\end{lemma}

\begin{proof}
By Definition \ref{def:HardyChar}, it suffices to find
a constant $C_\Gamma>0$ such that  inequality \eqref{eq:HardyChar} holds for every non-negative bounded $g\in LC(X)$ 
and every $x\in\Omega$ --- the remaining constants $\nu$ and $\kappa$ are given in the assumptions of the present lemma.
Fix such a function $g$ and a point $x\in\Omega$. Since
$g$ is bounded and inequality \eqref{eq:HardyChar} is invariant under multiplication
of $g$ with a strictly positive constant, we may further assume that
$g(y)\in [0,1]$ for all $y\in X$.
Then the desired estimate \eqref{eq:HardyChar},
with $C_\Gamma=C_\alpha$, follows immediately 
from Lemma \ref{lem:rewrite} and the assumptions.
\end{proof}

The converse of Lemma \ref{lem:alphaClass} is also true, as we will
see in Section \ref{s.self_improvement}.
Therein the following inequalities for the $\alpha$-function become useful.

\begin{lemma}\label{lem:alphaEst} 
Let  $\Omega\subsetneq X$  be an open set.  
Let  $0\le \tau<\tau'$, $\kappa>1$, $1\le p<\infty$ and $\nu>C_{\mathrm{QC}}$.  
Then
\[
\alpha_{p,w,\Omega}(\nu,\kappa,\tau) \leq \alpha_{p,w,\Omega}(\nu,\kappa,\tau'),\qquad 
\alpha_{p,w,\Omega}(\nu,\kappa,\tau) \leq \nu,
\]
and, for   every   $M\geq 1$, 
\[\alpha_{p,w,\Omega}(\nu,\kappa,M\tau) \leq M\alpha_{p,w,\Omega}(\nu,\kappa, \tau).
\]
\end{lemma}

\begin{proof}
These inequalities are clear from the definition 
of $ \alpha_{p,w,\Omega}(\nu,\kappa,\tau)$ in \eqref{eq:alphaDef}. The second inequality 
also uses the fact that $g$ is bounded by $1$
and quasiconvexity, that is, existence of a curve with $\len(\gamma)\le \nu d(x,\Omega^c)$.  
\end{proof}

\section{Self-improvement property for $p$-Hardy weights}\label{s.self_improvement}

In this section we examine self-improvement properties of
$p$-Hardy weights for $1<p<\infty$.
We assume that $w$ is a $p_0$-Poincar\'e weight
for some $p_0<p$. 
This assumption allows us to focus on the new phenomena
that arise especially in connection with the improvement of pointwise $p$-Hardy inequalities.
Recall that if the metric space $X$ is complete and $X$ supports a $(1,p)$-Poincar\'e
inequality, that is, \eqref{e.familiar} with $w=1$ holds for all balls $B\subset X$ whenever 
$u\in \Lip(X)$ and $g$ is an upper gradient of $u$, then
there exists $p_0<p$ such that $X$ supports a $(1,p_0)$-Poincar\'e
inequality; see~\cite{MR2415381}  and 
see also Lemma~\ref{lem.w_is_poinc} concerning this assumption
for distance weights in $\R^n$.
It is plausible that also $p$-Poincar\'e weights enjoy
self-improvement properties, but in the present work we
will not focus on this aspect. 

The following
Theorem \ref{t.converse_and_impro} implies
a self-improvement property for $p$-Hardy weights.
This result also provides a
converse of Lemma \ref{lem:alphaClass} for $p>1$.

\begin{theorem}\label{t.converse_and_impro}
Let $1<p_0<p<\infty$, let  $\Omega\subsetneq X$  be an open set and 
let $w$ be a semilocally doubling
weight for $\Omega$. Assume that $w$ is a 
$p_0$-Poincar\'e weight for $\Omega$
and a $p$-Hardy weight for $\Omega$.
Then there exist an exponent $q\in (p_0,p)$ and constants 
$N>C_{\mathrm{QC}}$, $K> 1$ and $C_\alpha>0$ such that
\begin{equation}\label{e.main_dest}
\alpha_{q,w,\Omega}(  N,K, \tau)\le C_\alpha\tau
\end{equation}
whenever $\tau\ge 0$. 
\end{theorem}

\begin{proof}
First, we fix some constants to give accurate bounds.
In Definition \ref{def:HardyChar}, inequality~\eqref{eq:HardyChar} holds with constants   $C_\Gamma>0$, $\nu_\Gamma>C_{\mathrm{QC}}$ and $\kappa_\Gamma> 1$. Also, denote by $C_{\mathrm{A}}>0$, $\nu_{\mathrm{A}} >C_{\mathrm{QC}}$ and $\kappa_{\mathrm{A}}> 1$ the constants from inequality \eqref{eq:PIchar} in Definition~\ref{def:PIchar}, for the exponent $p_0<p$. 
By H\"older's inequality we may assume $p/2<p_0$.
Without loss of generality, we may also assume that 
$\kappa_\Gamma=\kappa_{\mathrm{A}}=:\kappa$
and $\nu_\Gamma=\nu_{\mathrm{A}}=:\nu$.  

\vskip1mm
\textbf{Step 1: Estimate to prove, strategy and parameters.}
Assume that we have found parameters $k \in \N$, 
$K,S\in (1,\infty)$, $N \in (C_{\mathrm{QC}},\infty)$, $M>1$ and $\delta \in (0,1)$ such that,
for each $q\in (p_0,p)$ and every $\tau>0$, we have  
\begin{equation}\label{eq:iterationEst}
  \alpha_{q,w,\Omega}(N,  K ,\tau) \leq S \tau + \delta \max_{i=1, \dots, k} \big(M^{-iq/p} \alpha_{q,w,\Omega}(N,K,M^i \tau)\big).
\end{equation} 
From this inequality and Lemma \ref{lem:alphaEst}, we obtain
\[\alpha_{q,w,\Omega}(N,K,\tau) \leq S \tau + \delta M^{k\frac{p-q}{p}} \alpha_{q,w,\Omega}(N,K,\tau)\qquad   \text{ for all }q\in (p_0,p)\text{ and }\tau>0.
\]
Observe that the last term on the right is finite by Lemma \ref{lem:alphaEst}.
In order to absorb this term to the left-hand side, we need $\delta M^{k\frac{p-q}{p}} < 1$. 
This can be ensured  by  choosing  $q\in (p_0,p)$  so close to $p$ that 
\[0<p-q < \frac{p\ln(\frac{1}{\delta})}{k\ln(M)}.\]
With this choice of $q$ we find for all $\tau>0$ that
\[\alpha_{q,w,\Omega}(N,K,\tau) \leq \left(\frac{S}{1-\delta M^{k\frac{p-q}{p}}}\right)\tau=:C_\alpha \tau.\]
This inequality holds also for $\tau=0$, which is seen by using monotonicity property of the $\alpha$-function, see Lemma \ref{lem:alphaEst}. 
Thus, the desired inequality \eqref{e.main_dest} follows from \eqref{eq:iterationEst}.
Hence, it suffices to find parameters, as above, for which inequality \eqref{eq:iterationEst}
holds for every $q\in (p_0,p)$ and $\tau>0$.
 
We begin by fixing the auxiliary parameters
\[
K=2\kappa,\quad N=3\nu,\quad M = 4,\quad  \delta =  \frac{1}{6}. 
\]
We also choose $k\in\N$ so large that $  C_{\Gamma}^p  \frac{2^p {D(w,10\kappa)^4}}{k^{p-1}}<(\frac{\delta}{3\kappa})^p$.
The last parameter is given by
$S= 1+M^{k}  \nu +  3 C_{\mathrm{A}} M^{k}$. 
For what follows
$q\in (p_0,p)$ and $\tau>0$ are arbitrary.
 
Now, the overall strategy is to
construct, for any $x\in\Omega$ and any $g \in \mathcal{E}^{K,\tau}_{q,w,x,\Omega}$, a curve $\gamma \in \Gamma(X)^N_{x,\Omega^c}$ 
such that, for some $i_0  = 1, \dots, k$
\begin{equation}\label{e.desired_curve}
\int_\gamma g \,ds \leq S \tau d(x, \Omega^c) + \delta M^{-i_0 q/p} \alpha_{q,w,\Omega}(N,K,M^{i_0} \tau) d(x, \Omega^c).
\end{equation}
Estimating the right-hand side 
by the maximum over possible $i_0$, then 
dividing both sides 
by $d(x,\Omega^c)$, and finally
taking the supremum over $x$ and $g$ as above,  proves 
inequality  \eqref{eq:iterationEst}.
This strategy first involves choosing a good level $i_0$ along with some proto-curve $\gamma_0$ 
having a small integral, and then 
adjusting the curve at the level $i_0$ by filling in certain gaps.

\vskip1mm
\textbf{Step 2: Choosing a good level $i_0$ and the proto-curve $\gamma_0$.}
Fix $x\in\Omega$ and $g \in \mathcal{E}^{K,\tau}_{q,w,x,\Omega}$. For each $i\ge 1$, we write 
\begin{align*}
E_i &:= \{y \in \Omega  \mid  \M_{q,w,K}^{\kappa d(x,\Omega^c)} g(y) > M^i \tau\},
\end{align*} and define a bounded function $h\colon X\to [0,\infty)$ by setting 
\[
h = \frac{1}{k}\sum_{i=1}^k \mathbf{1}_{E_i} M^{iq/p}.
\]
Since $E_j \supset E_i$ if $j\le i$ and $p/2<p_0<q<p$, it follows that
\begin{align*}
h^p &\leq \frac{1}{k^p} \sum_{j=1}^k \bigg(\sum_{i=1}^j M^{iq/p}\bigg)^p \mathbf{1}_{E_j} \le \frac{2^p}{k^p} \sum_{j=1}^k \mathbf{1}_{E_j} M^{jq}.
\end{align*}
 In the final estimate, we also use the  choice $M=4$ to obtain the factor $2^p$. 
Observe that  $\mathbf{1}_{E_i}\in LC(X)$ 
since $E_i$ is open, for each $i=1,\ldots,k$, by the lower
semicontinuity of $\M_{q,w,K}^{\kappa d(x,\Omega^c)} g$. Hence, we have
$h\in  LC(X) $. 
By sublinearity and monotonicity of the maximal function,  Lemma \ref{lem:HLmaxmaxest}, and 
the assumption that $g \in \mathcal{E}^{K,\tau}_{q,w,x,\Omega}$, where $K=2\kappa$, we obtain
\begin{equation}\label{e.essential}
\begin{split}
 \left(\M_{p,w,\kappa} h (x) \right)^p &  \leq \frac{2^p}{k^p}\sum_{j=1}^k (\M_{1,w,\kappa} \mathbf{1}_{E_j}(x)) M^{jq} \\
	       &  \leq \frac{2^p }{k^p}\sum_{j=1}^k \frac{D(w,10\kappa)^4}{M^{jq}}M^{jq}\leq \frac{2^p{D(w,10\kappa)^4}}{k^{p-1}}.
\end{split}			       
\end{equation}
Then, by the choice of $k$ and estimate \eqref{e.essential}, we have
 \[
   C_\Gamma d(x,\Omega^c)\M_{p,w,\kappa} h (x) < \frac{\delta}{3\kappa}d(x,\Omega^c).
 \]
Therefore by Definition \ref{def:HardyChar}, with exponent $p$, there is a
curve $\gamma_0 \in \Gamma(X)^{\nu}_{x,\Omega^c}$, which is parametrized by the arc length
and defined on the interval $[0,\len(\gamma_0)]$, such that
\begin{equation}\label{e.h_app}
\frac{1}{k}\sum_{i=1}^k M^{iq/p}\int_{\gamma_0} \mathbf{1}_{E_i} \,ds = \int_{\gamma_0} h\,ds \leq \frac{\delta}{3\kappa} d(x,\Omega^c) 
\end{equation}
and
\begin{equation}\label{e.gamma0_Lestimate}
\len(\gamma_0) \leq \nu d(x,\Omega^c).
\end{equation}
Without loss of generality, we may assume that $\gamma_0([0,\len(\gamma_0))\subset \Omega$. 
By inequality \eqref{e.h_app}, there  exists $i_0 \in \{1,\ldots,k\}$  such that
\begin{equation}\label{e.apriori}
\int_{\gamma_0} \mathbf{1}_{E_{i_0}}\,ds \leq \frac{\delta}{3\kappa} M^{-i_0q/p} d(x,\Omega^c).
\end{equation}

\vskip1mm
\textbf{Step 3: Adjusting the curve at level ${i_0}$ by filling in gaps.}
Recall that the proto-curve $\gamma_0$ is parametrized by arc length.
Let $O = \gamma_0^{-1}(E_{i_0})$ and write $T= [0, \len(\gamma_0)] \setminus O$.
By the lower semicontinuity of $g$ and the definition of $E_{i_0}$ we have, for all $t \in T\setminus \{\len(\gamma_0)\}$,
\begin{equation}\label{e.K_estimate}
g(\gamma_0(t))\le \M_{q,w,K}^{\kappa d(x,\Omega^c)} g(\gamma_0(t))\le  M^{i_0} \tau. 
\end{equation} 
Since $E_{i_0}$ is open in $X$, the set $O$ is relatively open
in $[0,\len(\gamma_0)]$. Observe that $0\not\in O$ since $g\in \mathcal{E}^{K,\tau}_{q,w,x,\Omega}$. 
Likewise $\len(\gamma_0)\not\in O$ since $\gamma_0(\len(\gamma_0))\in \Omega^c$.   
Hence, we can write $O$ as a union of so-called gaps:
\begin{equation}\label{e.cases}
O = \bigcup_{i \in I} (a_i,b_i),
\end{equation}
where $I\subset \N$ is a finite or infinite indexing set. We also
write $x_i\defeq\gamma_0(a_i)$, $y_i\defeq\gamma_0(b_i)$ and 
$d_i\defeq d(x_i,y_i)$  for each $i\in I$.
There are two cases to consider: either $d_i< d(x_i,\Omega^c)/(3\kappa)$ for all $i\in I$ or there exists
$i\in I$ such that $d_i\ge d(x_i,\Omega^c)/(3\kappa)$.
The latter also includes the case when
$y_i\in \Omega^c$ for some $i \in I$.
In both cases 
the gaps $(a_i,b_i)$ are pairwise disjoint and
$0\le a_i<b_i\le \len(\gamma_0)$ for each $i\in I$.
By inequality \eqref{e.apriori}, we have
\begin{equation}\label{e.d_i_estimates}
\sum_{i\in I} d_i \leq \sum_{i\in I} \len(\gamma_0\vert_{[a_i,b_i]})=
\sum_{i\in I} \int_{\gamma_0\vert_{[a_i,b_i]}} \mathbf{1}_{E_{i_0}}\,ds
\le \int_{\gamma_0} \mathbf{1}_{E_{i_0}}\,ds \le 
\frac{\delta}{3\kappa} M^{-i_0 q/p} d(x,\Omega^c).
\end{equation}

For each $i$ we next define a filling curve $\gamma_i \co [a_i,b_i] \to X$ connecting $\gamma_0(a_i)$ and $\gamma_0(b_i)$. 

\vskip1mm
\textbf{Case 1: We have $d(x_i,y_i)=d_i<d(x_i,\Omega^c)/(3\kappa)$ for all $i\in I$.}
Fix $i\in I$. If $d_i=0$, we define $\gamma_i(t) = \gamma_0(a_i)=\gamma_0(b_i)$ for each $t \in [a_i,b_i]$. 
From now on we 
assume that $d_i>0$ and proceed as follows. Observe that $\kappa<K$ and $x_i,y_i\in\Omega\setminus E_{i_0} $. This gives
\begin{equation}\label{e.M_q}
\M_{q,w,\kappa}^{\kappa d(x,\Omega^c)} g(x_i) \leq M^{i_0} \tau\quad \text{ and }\quad \M_{q,w,\kappa}^{\kappa d(x,\Omega^c)} g(y_i) \leq M^{i_0} \tau.
\end{equation}
We apply Definition \ref{def:PIchar} to the points $x_i$ and $y_i$. 
After a reparametrization, this yields a curve
$\gamma_i \co [a_i,b_i] \to X$ such that
$\gamma_i(a_i)= x_i$, $\gamma_i(b_i)= y_i$, 
\begin{equation}\label{eq:lenestgammai}
\len(\gamma_i) \leq \nu d(x_i, y_i) = \nu  d_i,
 \end{equation}
and, by using also H\"older's inequality and the fact that $p_0< q$,
\begin{equation}\label{e.i_curve}
\begin{split}
    &\int_{\gamma_i} g \,ds \leq   C_{\mathrm{A}}  d(x_i,y_i)\left( \M_{q,w,\kappa}^{\kappa d(x_i,y_i)} g (x_i) + \M_{q,w,\kappa}^{\kappa d(x_i,y_i)} g (y_i)\right) +
    \underbrace{C_{\mathrm{A}} M^{i_0} \tau d(x_i,y_i)}_{>0}.
        \end{split}
\end{equation} 
Here $\kappa d(x_i,y_i) \leq \kappa d(x,\Omega^c)$,
since by \eqref{e.d_i_estimates} we have
\[d(x_i,y_i)=d_i   \le \sum_{i\in I} d_i\leq d(x,\Omega^c).\]
This estimate together with \eqref{e.M_q} and \eqref{e.i_curve} gives
   \begin{equation}\label{eq:intestgammai}
   \begin{split}
 \int_{\gamma_i} g \,ds   & \leq
   3 C_{\mathrm{A}} M^{i_0} \tau d_i.
    \end{split}    
\end{equation}

We define a curve $\gamma \co [0, \len(\gamma_0)] \to X$ by setting $\gamma(t) = \gamma_0(t)$ if $t \in T $ and $\gamma(t) = \gamma_i(t)$ if $t \in (a_i,b_i)$ for some $i\in I$ that is uniquely determined by $t$. Then, by the length  estimates \eqref{e.gamma0_Lestimate} and \eqref{eq:lenestgammai},
followed by inequality \eqref{e.d_i_estimates}, we obtain
\begin{align*}
\len(\gamma) &\leq \len(\gamma_0) + \sum_{i\in I} \len(\gamma_i) \\
&\leq \nu d(x,\Omega^c) +  \nu \sum_{i\in I} d_i  \leq   2 \nu  d(x,\Omega^c) \leq N d(x,\Omega^c).
\end{align*}
From this it follows that $\gamma \in\Gamma(X)^N_{x,\Omega^c}$;  we remark that the required continuity and connecting
properties of $\gamma$ are straightforward to establish,
and we omit the details. Also, by 
inequalities \eqref{e.gamma0_Lestimate}, \eqref{e.K_estimate}, \eqref{e.d_i_estimates} and \eqref{eq:intestgammai},  we have
  \begin{align*}
      \int_\gamma g \,ds &
         =  \int_{T} g(\gamma_0(t)) \,dt  + 
      \sum_{i\in I} \int_{\gamma_i} g \,ds \\
                        &  \leq M^{i_0} \tau    \nu   d(x,\Omega^c) + 3 C_{\mathrm{A}} M^{i_0}  \tau  d(x,\Omega^c) \\
                        & \leq (M^{i_0}   \nu  +  3 C_{\mathrm{A}}  M^{i_0})\tau d(x,\Omega^c) 
                        \le S\tau d(x,\Omega^c).
 \end{align*}
Thus curve $\gamma$ satisfies inequality \eqref{e.desired_curve} and 
this concludes the proof of the first case.  
 
\vskip1mm
\textbf{Case 2: There exists $i\in I$ such that $d(x_i,y_i)= d_i\ge d(x_i,\Omega^c)/(3\kappa)$.} 
This includes the case when $b_i=\len(\gamma_0)$ for
some $i\in I$.
Write
\[
t=\inf\{a_i\mid i\in I\text{ and }d(x_i,y_i)\ge d(x_i,\Omega^c)/(3\kappa)\}\in [0,\len(\gamma_0)).
\]
The infimum is reached, that is, there exists an index
$i_0\in I$ such that $t=a_{i_0}$ and $d(x_{i_0},y_{i_0})\ge  d(x_{i_0},\Omega^c)/(3\kappa)$. 
Indeed, otherwise
there would exist a strictly decreasing sequence
$(a_{i_k})_{k\in \N}$ such that $i_k\in I$ and $d(x_{i_k},y_{i_k})\ge d(x_{i_k},\Omega^c)/(3\kappa)$
for all $k\in\N$, and $\lim_{k\to\infty} a_{i_k}=t$.
Clearly $a_{i_{k-1}}-a_{i_k}\to 0$ as $k\to\infty$. 
Since $\gamma_0$ is parametrized
by arc length, we obtain for all $k>1$
\begin{align*}
d(\gamma_0(a_{i_k}),\Omega^c)/(3\kappa)&=d(x_{i_k},\Omega^c)/(3\kappa)
\le d(x_{i_k},y_{i_k})\\
&=d(\gamma_0(a_{i_k}),\gamma_0(b_{i_k}))
\le b_{i_k}-a_{i_k}\le a_{i_{k-1}}-a_{i_k}\xrightarrow{k\to\infty} 0.
\end{align*}
Hence, by continuity, we have 
$d(\gamma_0(t),\Omega^c)=\lim_{k\to\infty} d(\gamma_0(a_{i_k}),\Omega^c)=0$.
Since $\Omega^c$ is closed, this implies
$\gamma_{0}(t)\in \Omega^c$. This is a contradiction, since
$t<\len(\gamma_0)$ and, on the other hand, 
we have assumed that $\gamma_0([0,\len(\gamma_0))\subset \Omega$.

Let $J\defeq\{i\in I\mid a_i<a_{i_0}\}$. Then
$d_i<d(x_i,\Omega^c)/(3\kappa)$ for all $i\in J$.
As in the previous case, for each $i\in J$, we can first construct curves $\gamma_i \co [a_i,b_i] \to X$ such that 
\begin{equation}\label{eq:lenestgammaiR}
   \len(\gamma_i) \leq  \nu d(x_i, y_i) = \nu   d_i,
\end{equation}
and
\begin{equation}\label{eq:intestgammaiR}
  \int_{\gamma_i} g \,ds \leq 3 C_{\mathrm{A}} M^{i_0} \tau d_i.
\end{equation}
For $i=i_0$ we are too close to the boundary and must proceed more carefully.
By using \eqref{e.d_i_estimates} and the 
equality $3K\delta= \kappa$, we first observe that
\[
Kd(x_{i_0},\Omega^c)\le 3\kappa Kd(x_{i_0},y_{i_0})
=3\kappa Kd_{i_0}\le  3K \delta d(x,\Omega^c)\le \kappa d(x,\Omega^c).
\]
We still have that $x_{i_0} \in \Omega\setminus E_{i_0} $, and thus
\[
\M_{q,w,K} g(x_{i_0}) \le \M_{q,w,K}^{\kappa d(x,\Omega^c)} g(x_{i_0})
\leq M^{i_0} \tau.
\] 
From this it follows that $g\in \mathcal{E}^{K,M^{i_0} \tau}_{q,w,x_{i_0},\Omega}$.
By definition \eqref{eq:alphaDef} of the function $\alpha_{q,w,\Omega}(N,K,M^{i_0} \tau)$, 
we obtain a curve $\gamma_{i_0}\co [a_{i_0},b_{i_0}]\to X$ connecting $x_{i_0} \in \Omega$ to $\Omega^c$ such that
\begin{equation}\label{eq:lenestgammaiInit}
\len(\gamma_{i_0}) \leq   Nd(x_{i_0},\Omega^c) \le 3\kappa N  d(x_{i_0},y_{i_0}) = 3\kappa N  d_{i_0}
\end{equation}
and  
\begin{equation}\label{eq:intestgammaiR_0}
\begin{split}
\int_{\gamma_{i_0}} g\,ds &\leq d(x_{i_0},\Omega^c)\alpha_{q,w,\Omega}(N,K,M^{i_0}\tau) +   \underbrace{\tau d(x,\Omega^c)}_{>0}  \\
&\le  3\kappa   d_{i_0} \alpha_{q,w,\Omega}(N,K,M^{i_0} \tau)+\tau d(x,\Omega^c).
\end{split}
\end{equation}

We now define a curve $\gamma \co [0, b_{i_0}] \to X$ by setting  $\gamma(t) = \gamma_0(t)$
if $t \in T\cap [0, a_{i_0}]$, 
$\gamma(t) = \gamma_i(t)$ if $t \in (a_i,b_i)$ for some $i\in J$, which is uniquely determined by $t$, 
and $\gamma(t)= \gamma_{i_0}(t)$ for every $t \in (a_{i_0},b_{i_0}]$. Then by 
\eqref{e.gamma0_Lestimate},  \eqref{e.d_i_estimates}, \eqref{eq:lenestgammaiR}, \eqref{eq:lenestgammaiInit}, and 
our choices of $N$ and $\delta$, we obtain 
 \begin{align*}
 \len(  \gamma ) &\leq \len(\gamma_0) +   \sum_{i\in J} \len(\gamma_i) +\len(\gamma_{i_0})\\
 &\leq (  \nu  +   \nu + \delta N) d(x,\Omega^c) \leq N d(x,\Omega^c),
 \end{align*}
and thus $\gamma\in \Gamma(X)^N_{x,\Omega^c}$. 
Finally, by inequalities  \eqref{e.gamma0_Lestimate}, \eqref{e.K_estimate}, \eqref{e.d_i_estimates}, \eqref{eq:intestgammaiR}, and \eqref{eq:intestgammaiR_0} we have 
\begin{align*}
\int_{ \gamma } g\,ds    &=  \int_{T\cap [0,a_{i_0}]} g(\gamma_0(t))\,dt   + 
\sum_{i\in J} \int_{\gamma_i} g \,ds + \int_{\gamma_{i_0}} g \,ds  \\
                    & \leq  M^{i_0} \tau    \nu   d(x,\Omega^c) + 
                     3 C_{\mathrm{A}} M^{i_0}  \tau   d(x,\Omega^c) + 3\kappa d_{i_0} \alpha_{q,w,\Omega}(N,K,M^{i_0} \tau)+  \tau d(x,\Omega^c)  \\
                        & \leq  S\tau d(x,\Omega^c) + \delta M^{-i_0 q/p} \alpha_{q,w,\Omega}(N,K,M^{i_0} \tau) d(x,\Omega^c). 
\end{align*}
This shows that \eqref{e.desired_curve} holds also in Case 2,
and the proof is complete.
\end{proof}

\section{Pointwise $(p,w)$-Hardy inequalities}\label{s.pw_hardy}

The definition of a pointwise $(p,w)$-Hardy inequality is as follows; recall that $\Omega^c=X\setminus \Omega$.

\begin{definition}\label{def:ptwisePHardy}
Let $1\le p<\infty$, let  $\Omega\subsetneq X$ be an open set, and let $w$ be a weight
for $\Omega$. We say that a pointwise $(p,w)$-Hardy inequality holds in $\Omega$ if there 
exist constants $C_{\mathrm{H}}>0$ and  $\kappa>1$ such that for every Lipschitz function $u \in \Lip_{0}(\Omega)$, 
every bounded upper gradient $g$ of $u$ and every $x \in \Omega$, we have
\begin{equation}\label{e.pointwise}
\lvert u(x)\rvert \leq  C_{\mathrm{H}}\, d(x,\Omega^c) \M_{p,w,\kappa} g (x).
\end{equation}
\end{definition}

These pointwise $(p,w)$-Hardy inequalities in fact characterize the class of $p$-Hardy weights for $\Omega$,
thus explaining the terminology.

\begin{lemma}\label{lem:HardyChar} 
Let $1\le p<\infty$, let  $\Omega\subsetneq X$  be an open set and let $w$ be a semilocally doubling 
weight for $\Omega$.
Then a pointwise $(p,w)$-Hardy inequality holds in $\Omega$ if, and only if, 
$w$ is a $p$-Hardy weight for $\Omega$.
\end{lemma}

\begin{proof}
Throughout this proof, we tacitly assume that
curves are parametrized by arc length. 
First assume that a pointwise $(p,w)$-Hardy inequality \eqref{e.pointwise} holds in $\Omega$ with constants
$C_{\mathrm{H}}>0$ and $\kappa_\Gamma>1$.  
Let $g\in LC(X)$ be a non-negative and bounded function, 
and fix $x\in \Omega$ and $\delta>0$.
We define a function $u\colon X\to [0,\infty)$ by setting  
\begin{equation}\label{e.u_def}
  u(y) = \inf_{\gamma} \int_\gamma h \,ds,\qquad   y\in X,  
\end{equation}
where 
$h = g + \M_{p,w,\kappa_\Gamma} g (x)+\delta$
and the infimum is taken over all curves
$\gamma$ in $X$ connecting $y$ to $\Omega^c$;   
note that $h$ is a non-negative bounded Borel function. 
Clearly, we have $u=0$ in $\Omega^c$.
Fix $y,z\in X$ and consider any curve $\sigma$ connecting $y$ to $z$.
As in Step 2 of the proof of Lemma~\ref{l.poinc_suff},
we assume that $u(y)>u(z)$ and 
fix $\eps>0$. We let $\gamma$ be a curve in $X$ that
connects $z$ to $\Omega^c$ and satisfies inequality  
\[
u(z)\ge \int_\gamma h\,ds - \eps,
\]
and define $\sigma \gamma$ to be the concatenation of $\sigma$ and $\gamma$. Then,
as in the proof of Lemma \ref{l.poinc_suff},
\begin{align*}
\lvert u(y)-u(z)\rvert 
\le \int_{\sigma} h\,ds+\eps,
\end{align*}
and by taking $\eps\to 0_+$ we obtain
\begin{equation}\label{e.des}
\lvert u(y)-u(z)\rvert\le \int_\sigma h\,ds.
\end{equation}
This shows that $h$ is an upper gradient of $u$. 
Moreover, since $X$ is quasiconvex and $h$ is bounded, it follows from \eqref{e.des} that $u\in \Lip_0(\Omega)$.

Now, applying the assumed pointwise $(p,w)$-Hardy inequality \eqref{e.pointwise}
to $u\in \Lip_0(\Omega)$ and its bounded upper gradient $h$ yields
\[
u(x)
\le C_{\mathrm{H}}\,d(x,\Omega^c)\mathcal{M}_{p,w,\kappa_\Gamma} h(x)<\infty.
\]
Since $u(x)\ge \delta d(x,\Omega^c)>0$, by \eqref{e.u_def} there is a curve $\gamma$ in $X$ 
  connecting $x$ to $\Omega^c$ such that  
\begin{equation}\label{e.strong}
\begin{split}
\int_\gamma g\,ds + (\mathcal{M}_{p,w, \kappa_\Gamma} g(x)+\delta)\len(\gamma)
 &= \int_\gamma h\,ds   \le 2u(x)\\
&\le 2C_{\mathrm{H}}\,d(x,\Omega^c)(\mathcal{M}_{p,w,\kappa_\Gamma} h(x))
\\&\le 2C_{\mathrm{H}}\,d(x,\Omega^c)(2\mathcal{M}_{p,w,\kappa_\Gamma} g(x)+\delta)
\\&\le 4C_{\mathrm{H}}\,d(x,\Omega^c)(\mathcal{M}_{p,w,\kappa_\Gamma} g(x)+\delta).
\end{split}
\end{equation}
The last inequality follows from the sublinearity of maximal function. 
We can now conclude from \eqref{e.strong} that
$\len(\gamma)\le 4C_{\mathrm{H}}\,d(x,\Omega^c)$.
By taking $\delta\to 0_+$, we also obtain from \eqref{e.strong} that
\eqref{eq:HardyChar} holds, that is,
\[
\inf_{\gamma \in \Gamma(X)^\nu_{x,\Omega^c}} \int_{\gamma} g \,ds \leq C_\Gamma \,d(x,\Omega^c)  \M_{p,w,\kappa} g (x)
\]
with 
\[C_\Gamma=4C_{\mathrm{H}},\qquad \kappa=\kappa_\Gamma,\qquad \nu>\max\{C_{\mathrm{QC}},4C_{\mathrm{H}}\}.\]

For the converse implication, we assume that inequality \eqref{eq:HardyChar} holds for all non-negative and bounded $g\in LC(X)$ and for all $x\in\Omega$. 
We need to prove that a pointwise $(p,w)$-Hardy inequality holds in $\Omega$.
To this end, we fix $x\in \Omega$ and $u\in \Lip_0(\Omega)$, and let $g$
be a bounded upper gradient of $u$.
Since $g$ is not necessarily lower semicontinuous, 
some approximation is first needed so that we 
get to apply \eqref{eq:HardyChar} and thereby establish inequality \eqref{e.pointwise}.

Let
$(g_N)_{N\in\N}$ be a pointwisely increasing sequence of simple Borel functions such that 
$\lim_{N\to\infty} g_N=g$ uniformly in $X$. Fix $\varepsilon>0$. 
By the uniform convergence, there exists $N\in\N$ such that
for all $\gamma\in \Gamma(X)^\nu_{x,\Omega^c}$ we have
\begin{equation}\label{e.simple}
\begin{split}
 \int_{\gamma} g \,ds &=\, \int_{\gamma} g_N \,ds + \int_{\gamma} (g-g_N) \,ds \\
 &\le 
 \int_{\gamma} g_N \,ds + \sup_{y\in X}(g(y)-g_N(y))\len(\gamma)\\
  &\le 
 \int_{\gamma} g_N \,ds + \sup_{y\in X}(g(y)-g_N(y))\nu d(x,\Omega^c)\\
 &\le  \int_{\gamma} g_N \,ds+\varepsilon.
 \end{split}
\end{equation}
Let $g_{N,x,\varepsilon}\in LC(X)$ be 
the non-negative bounded approximant of $g_N$ given by Lemma \ref{lem:approx} with
$F=\{x\}$.
By inequality \eqref{eq:HardyChar} and Lemma \ref{lem:approx},
there exists $\gamma_{N}\in \Gamma(X)^\nu_{x,\Omega^c}$ 
defined on the interval $[0,\ell(\gamma_N)]$ such that
\begin{equation}\label{e.addeneum}
\begin{split}
\int_{\gamma_{N}} g_{N,x,\varepsilon} \,ds&\leq C_\Gamma d(x,\Omega^c) \M_{p,w,\kappa} g_{N,x,\varepsilon} (x) +\varepsilon \\&\leq  C_\Gamma d(x,\Omega^c) \left( \M_{p,w,\kappa} g_{N} (x)+\varepsilon \right)+\varepsilon\\&\leq C_\Gamma d(x,\Omega^c) \left( \M_{p,w,\kappa} g (x)+\varepsilon \right)+\varepsilon.
\end{split}
\end{equation}
 
Without loss of generality, we may assume that $\gamma_N(t)=x$ only if $t=0$. 
On the other hand, by Lemma \ref{lem:approx}, we have  $g_{N}\le g_{N,x,\varepsilon}$ in $X\setminus \{x\}$. Inequalities \eqref{e.simple} and \eqref{e.addeneum}, with $\gamma=\gamma_N$, imply that
\begin{align*}
\int_{\gamma_N} g \,ds &\leq   \int_{\gamma_N} g_N \,ds + \varepsilon 
\le  \int_{\gamma_N} g_{N,x,\varepsilon} \,ds + \varepsilon \\
&\le  C_\Gamma d(x,\Omega^c) \left( \M_{p,w,\kappa} g (x)+\varepsilon \right)+2\varepsilon. 
\end{align*}
Since $g$ is an upper gradient of $u\in \Lip_0(\Omega)$ and $\gamma_{N}(\len(\gamma_{N}))\in \Omega^c$, 
we obtain
\begin{align*}
 \lvert u(x)\rvert 
 &=\lvert u(\gamma_N(0))-u(\gamma_{N}(\len(\gamma_{N})))\rvert\\&
 \leq \int_{\gamma_N} g \,ds\le C_\Gamma d(x,\Omega^c) \left( \M_{p,w,\kappa} g (x)+\varepsilon \right)+2\varepsilon,
\end{align*} 
and letting $\varepsilon \to 0_+$ gives the pointwise $(p,w)$-Hardy inequality \eqref{e.pointwise}
with $C_{\mathrm{H}}=C_\Gamma$ and $\kappa$.  
\end{proof}

\begin{remark}
Let  $\Omega\subsetneq X$  be an open set and let $w$ be 
a weight for $\Omega$ such
that a pointwise $(p,w)$-Hardy inequality holds
in $\Omega$, with constants $C_{\mathrm{H}}>0$ and  $\kappa  > 1$.
Then the proof of Lemma \ref{lem:HardyChar}, with $g=0$, shows that for every 
$\eps>0$ and every $x\in\Omega$ there exists a curve $\gamma$ 
that connects $x$ to $\Omega^c$ in $X$ such that
$\len(\gamma)\le (1+\eps)C_{\mathrm{H}}d(x,\Omega^c)$.
\end{remark}

The following is our main result.

\begin{theorem}\label{t.main_improvement}
Let $1\le p_0<p<\infty$, let  $\Omega\subsetneq X$  be an open set,
and assume that $w$ is a semilocally doubling $p_0$-Poincar\'e weight
for $\Omega$.
If a pointwise $(p,w)$-Hardy inequality holds in $\Omega$,
then there exists $q\in (p_0,p)$ such that  
a pointwise $(q,w)$-Hardy inequality holds in $\Omega$. 
\end{theorem}

\begin{proof}
By Lemma \ref{lem:HardyChar}, we find that
$w$ is a $p$-Hardy weight for $\Omega$. Theorem~\ref{t.converse_and_impro} and Lemma~\ref{lem:alphaClass}
imply that there exists $q\in (p_0,p)$ such that  
$w$ is a $q$-Hardy weight for $\Omega$.
Lemma~\ref{lem:HardyChar} implies
that a pointwise $(q,w)$-Hardy inequality holds in $\Omega$.
\end{proof}

\begin{remark}
The proofs of the results in Sections~\ref{s.weights}, \ref{s.poinc_weight}, \ref{s.self_improvement} and~\ref{s.pw_hardy}
show that the semilocal doubling property in Definition~\ref{def.semilocal} is not
really needed to hold for every $\kappa>0$ but for every $0<\kappa\le\kappa_0$ 
with a large enough $\kappa_0$ depending on
the parameters in the assumed $p_0$-Poincar\'e weight property and
pointwise $(p,w)$-Hardy inequality.
\end{remark}

\section{Applications}\label{s.applications}

In this section we show how the self-improvement of pointwise $(p,w)$-Hardy
inequalities can be applied in the context of integral versions of
weighted Hardy inequalities. Here we need to know, for all $1<q<\infty$ and all $0<\kappa<\infty$, the
boundedness of the restricted weighted maximal operator 
$\M_{1,w,\kappa}\colon L^q(X;w\,d\mu) \to L^q(\Omega;w\,d\mu)$, where $w$ is a semilocally
doubling weight for an open set $\Omega$. 
If $w$ is a doubling weight in $X$, then this $L^q$-boundedness of $\M_{1,w,\kappa}$
follows from the maximal function theorem in $X$; see, for instance~\cite[Section~3.2]{MR2867756}.
In our case the weight $w$ is not necessarily doubling, but the boundedness follows with
a suitable adaptation of the proof of the doubling case, given by the following lemma.

\begin{lemma}\label{l.maximal_bounded}
Let $0<\kappa<\infty$ and $1< q<\infty$, let $\Omega\subsetneq X$  be an open set,
and assume that $w$ is a semilocally doubling weight
for $\Omega$. Then the restricted weighted maximal operator 
${\M}_{1,w,\kappa}\colon L^q(X;w\,d\mu) \to L^q(\Omega;w\,d\mu)$ is
bounded, that is, there is a constant $C=C_{q,\kappa,w}$ such that
\[
\int_{\Omega} \bigl(\M_{1,w,\kappa} f\bigr)^{q} w\,d\mu
\le C\int_{X} \lvert f\rvert^q w\,d\mu,
\]
for every $f\in L^q(X;w\,d\mu)$. 
\end{lemma}

\begin{proof} Clearly $\M_{1,w,\kappa} \co L^\infty(X;w\,d\mu) \to L^\infty(\Omega;w\,d\mu)$ is bounded.
Hence by interpolation it suffices to prove that $\M_{1,w,\kappa}$ is of corresponding weak type $(1,1)$, 
compare to the proof of~\cite[Theorem~3.13]{MR2867756}. Let $f \in L^1(X;w\,d\mu)$ and $0<\tau<\infty$. 
We estimate the $w\,d\mu$-measure of $E=\{x\in \Omega \ |\  \M_{1,w,\kappa} f(x) > \tau \}$. The set $E$ has a cover by balls in
\[
\mathcal{B}=\biggl\{ B=B(x,r) \ \Big|\  x \in \Omega, \  \frac{1}{w(B)} \int_B |f|w\,d\mu > \tau ,\  0<r< \kappa d(x,\Omega^c)\biggr\}.
\]
By the $5r$-covering lemma \cite[Lemma 1.7]{MR2867756}, we obtain a countable subfamily 
$\mathcal{B}' \subset \mathcal{B}$ of pairwise disjoint balls such that
\[
E \subset \bigcup_{B \in \mathcal{B'}} 5B.
\]
Then
\begin{align*}
w(E) &\leq \sum_{B(x,r) \in \mathcal{B'}} w(B(x,5r)) \leq D(w,5\kappa)^3 \sum_{B(x,r) \in \mathcal{B'}}w(B(x,r)) \\
&\leq D(w,5\kappa)^3\sum_{B\in \mathcal{B'}} \frac{\int_B |f|w\,d\mu}{\tau}\leq 
\frac{D(w,5\kappa)^3}{\tau}\int_X \lvert f\rvert w\,d\mu.
\end{align*}
On the first line we used semilocal doubling with $x\in\Omega$ and $0<r<\kappa d(x,\Omega^c)$, 
and on the last line we used the fact that the balls in $\mathcal{B}'$ are pairwise disjoint.
This shows the desired weak type $(1,1)$ property and the proof is completed by interpolation.
\end{proof}

\begin{theorem}\label{t.abstract}
Let $1\le p_0<p<\infty$, let $\Omega\subsetneq X$ be an open set,
and assume that $w$ is a semilocally doubling $p_0$-Poincar\'e weight
for $\Omega$.
Assume that a pointwise $(p,w)$-Hardy inequality holds in $\Omega$.
Then there exists a constant $C>0$ such that
the $(p,w)$-Hardy inequality
\[
\int_\Omega\frac{\lvert u(x)\rvert^p}{d(x,\Omega^c)^p}\,w(x)\,d\mu(x)
\le C\int_\Omega g(x)^p w(x)\,d\mu(x)
\]
holds for every $u\in \Lip_0(\Omega)$ and for all bounded upper
gradients $g$ of $u$.
\end{theorem}

\begin{proof}
By Theorem~\ref{t.main_improvement} there exists $q\in (p_0,p)$
such that a pointwise $(q,w)$-Hardy inequality holds in $\Omega$ with $1<\kappa<\infty$.
Let $u\in \Lip_0(\Omega)$ and let $g$ be a bounded upper
gradient of $u$.
Without loss of generality, we may assume that $g=0$ in $\Omega^c$.
It is immediate that
\begin{equation*}
\bigl({\M}_{q,w,\kappa} g(x)\bigr)^p
= \bigl({\M}_{1,w, \kappa} g^q(x)\bigr)^{\frac p q},
\end{equation*}
for every $x\in\Omega$,
and on the other hand the pointwise $(q,w)$-Hardy inequality, raised to power $p$,
implies
\[
\frac{\lvert u(x)\rvert^p}{d(x,\Omega^c)^p}\le  C\bigl(\M_{q,w,\kappa}g(x)\bigr)^p
\]
for every $x\in\Omega$.
Since $p/q>1$, by the $L^{p/q}$-boundedness of $\M_{1,w,\kappa}$ from Lemma~\ref{l.maximal_bounded} we obtain
\begin{align*}
\int_\Omega \frac{\lvert u(x)\rvert^p}{d(x,\Omega^c)^p}\, w(x)\,d\mu(x)
& \le C\int_\Omega \bigl(\M_{q,w,\kappa}g(x)\bigr)^p w(x)\,d\mu(x)\\
& = C\int_{\Omega} \bigl(\M_{1,w,\kappa} g^q (x)\bigr)^{\frac p q} w(x)\,d\mu(x)\\
&\le C\int_{X} g(x)^p w(x)\,d\mu(x)
= C\int_\Omega g(x)^p w(x)\,d\mu(x),
\end{align*}
and this proves the claim. 
\end{proof}

Next we concentrate on the special case where $X=\R^n$ is equipped with the Euclidean distance and the Lebesgue measure,
and weights $w$ are powers of the distance function $x\mapsto d(x,\Omega^c)$.
Let $1\le p<\infty$ and let  
$\Omega\subsetneq\R^n$ be an open set.
We say that a pointwise $p$-Hardy inequality
holds in $\Omega$, if there exists a constant $C>0$ such that
\begin{equation}\label{e.pw_hardy}
\lvert u(x)\rvert\le C d(x,\Omega^c)
\bigl(\M_{2d(x,\Omega^c)}\lvert\nabla u\rvert^p(x)\bigr)^{\frac1p},
\end{equation}
for every $x\in\Omega$ and every $u\in \Lip_0(\Omega)$.
Here $\M_{2d(x,\Omega^c)}$ is the usual restricted maximal operator, which corresponds to 
$\M_{1,1,2}$ in the notation introduced in Section~\ref{s.weights}. These pointwise inequalities were introduced and studied by
Haj\l asz~\cite{MR1458875} and Kinnunen and Martio~\cite{MR1470421}, 
and they can be regarded as pointwise variants of the usual $p$-Hardy inequality
\begin{equation}\label{e.hardy}
\int_\Omega\frac{\lvert u(x)\rvert^p}{d(x,\Omega^c)^p}\,dx
\le C\int_\Omega\lvert \nabla u(x)\rvert^p\,dx.
\end{equation} 

If~\eqref{e.pw_hardy} holds for a function
$u\in \Lip_0(\Omega)$ at every $x\in\Omega$, but with an exponent $1< q<p$, then
the maximal function theorem implies that
\eqref{e.hardy} with exponent $p$ holds for $u$ with a constant $C$ independent of $u$.
However, the passage from~\eqref{e.pw_hardy} to \eqref{e.hardy}, with the same exponent $1<p<\infty$,
is not at all obvious. This was established in~\cite{MR2854110} using an indirect route, first
showing the equivalence between the validity of~\eqref{e.pw_hardy} and the uniform $p$-fatness
of $\Omega^c$, and then applying the known self-improvement of the latter, which in $\R^n$ is
by Lewis~\cite{MR946438} and in metric spaces by Bj\"orn, MacManus and Shanmugalingam~\cite{MR1869615}. A direct proof for the self-improvement of
pointwise $p$-Hardy inequalities, which applies also in metric spaces, was recently given in~\cite{MR3976590}. 

The following weighted version of the pointwise $p$-Hardy inequality was
considered in \cite{KoskelaLehrback2009}: 
\begin{equation}\label{e.pw_beta_hardy}
\lvert u(x)\rvert\le C d(x,\Omega^c)^{1-\frac{\beta}{p}}
\bigl(\M_{2d(x,\Omega^c)}\bigl(\lvert\nabla u\rvert^q d(\cdot,\Omega^c)^{\frac{\beta q}{p}}\bigr)(x)\bigr)^{\frac1q},
\end{equation}
for every $x\in\Omega$ and every $u\in \Lip_0(\Omega)$, where $1<q<p$ are fixed.
As in the unweighted case, with an application of the maximal function theorem for exponent $\frac p q > 1$,
this implies the weighted $(p,\beta)$-Hardy inequality 
\begin{equation}\label{e.beta_hardy}
\int_\Omega \lvert u(x)\rvert^p d(x,\Omega^c)^{\beta-p}\,dx
\le C\int_\Omega\lvert \nabla u(x)\rvert^p  d(x,\Omega^c)^{\beta}\,dx.
\end{equation} 
A more natural formulation for the weighted pointwise Hardy inequality~\eqref{e.pw_beta_hardy} 
would  have been with $q=p$, 
but then the passage to inequality~\eqref{e.beta_hardy} would not have been possible
with a direct use of the maximal function theorem. 

Now, using the general technology developed in this paper, we can show that 
the validity of \eqref{e.pw_beta_hardy},
with $1<q=p<\infty$, implies \eqref{e.beta_hardy}, at least in the case $\beta\ge 0$. 
We begin by proving that 
in this case the weight $w(x)=d(x,\Omega^c)^\beta$, for $x\in\R^n$,
is a semilocally doubling $p_0$-Poincar\'e weight for $\Omega$, for every $1\le p_0<\infty$.

\begin{lemma}\label{lem.w_is_poinc}
Let $1\le p_0<\infty$ and $\beta\ge 0$, and let $\Omega\subsetneq\R^n$ be an open set.
Define $w(x)=d(x,\Omega^c)^\beta$ for all $x\in\R^n$. 
Then $w$ is a semilocally doubling $p_0$-Poincar\'e weight for $\Omega$.
\end{lemma}

\begin{proof}
Let $\kappa> 0$, $x\in\Omega$ and $0<r\le \kappa d(x,\Omega^c)$.  
There exists $C=C(n,\beta,\kappa)$ such that
\begin{equation}\label{e.equiv_moved}
\begin{split}
C^{-1}r^n d(x,\Omega^c)^\beta &\le  w(B(x,r))
=\int_{B(x,r)} d(y,\Omega^c)^\beta\,dy\le Cr^nd(x,\Omega^c)^\beta,
\end{split}
\end{equation}
and this shows that $w$ is a semilocally doubling weight for $\Omega$.

To prove the $p_0$-Poincar\'e weight property, we 
let $u\in \Lip(\R^n)$.
There exists a bounded upper gradient $g_u$ of $u$ such that
\begin{equation}\label{e.minimality}
g_u=\lvert \nabla u\rvert\qquad \text{ and }\qquad \lvert \nabla u\rvert\le g
\end{equation} almost everywhere in $\R^n$ whenever $g$ is a bounded
upper gradient of $u$; 
we refer to the proof of~\cite[Corollary~1.47]{MR2867756} and~\cite[Proposition~A.3]{MR2867756}.
Let $x\in\Omega$ and let $B=B(x,r)$ be a ball
with $2B=B(x,2r)\subset \Omega$.
We have $0<r\le d(x,\Omega^c)/2$ and $d(y,\Omega^c)\le 2d(x,\Omega^c)\le 4d(y,\Omega^c)$ for every $y\in B$. 
By \eqref{e.equiv_moved}, with $\kappa=1/2$, and the well-known  $1$-Poincar\'e inequality in $\R^n$, 
we have 
\begin{align*}
&\frac{1}{w(B)}\int_{B} \lvert u(y)-u_{B;w}\rvert w(y) \,dy\le \frac{2}{w(B)}\int_{B} \lvert u(y)-u_{B;1}\rvert w(y) \,dy\\
&\qquad \le C(\beta)\frac{d(x,\Omega^c)^\beta}{w(B)}\int_{B} \lvert u(y)-u_{B;1}\rvert \,dy
\le \frac{C(n,\beta)}{\lvert B\rvert}\int_{B}\lvert u(y)-u_{B;1}\rvert \,dy\\
&\qquad \le C(n,\beta)\frac{r}{\lvert B\rvert} \int_B \lvert \nabla u(y)\rvert\,dy
\le C(n,\beta)\frac{r}{w(B)} \int_B \lvert \nabla u(y)\rvert w(y) \,dy \\
&\qquad \le C(n,\beta)\frac{r}{w(B)} \int_B g(y) w(y)\,dy
\end{align*}
whenever $g$ is a bounded upper gradient of $u$,
where the final step follows from the second inequality in \eqref{e.minimality}.
This together with H\"older's inequality and Lemma \ref{l.poinc_suff}, 
with $\lambda=1$, proves that $w$ is a $p_0$-Poincar\'e weight in $\Omega$.
\end{proof}

The claim that weighted pointwise $(p,\beta)$-Hardy inequality~\eqref{e.pw_beta_hardy_natural}, with $\beta\ge 0$,
implies the integral version of the $(p,\beta)$-Hardy inequality is now a special case of Theorem~\ref{t.abstract}.

\begin{theorem}\label{t.euclidean}
Let $1<p<\infty$ and $\beta \ge 0$, and
let $\Omega\subsetneq \R^n$ be an open set.
Assume that there exists a constant $C>0$ such that
\begin{equation}\label{e.pw_beta_hardy_natural}
\lvert u(x)\rvert\le C d(x,\Omega^c)^{1-\frac{\beta}{p}}
\bigl(\M_{2d(x,\Omega^c)}\bigl(\lvert\nabla u\rvert^p d(\cdot,\Omega^c)^{\beta}\bigr)(x)\bigr)^{\frac1p},
\end{equation}
for every $x\in\Omega$ and every $u\in \Lip_0(\Omega)$.
Then the weighted $(p,\beta)$-Hardy inequality~\eqref{e.beta_hardy}
holds for every $u\in \Lip_0(\Omega)$, with a constant independent of $u$.
\end{theorem}

\begin{proof}
Define $w(x)=d(x,\Omega^c)^\beta$ for every $x\in\R^n$ and let $\kappa=2$.
By Lemma~\ref{lem.w_is_poinc}, $w$ is a semilocally doubling $1$-Poincar\'e weight for $\Omega$.
Let $u\in \Lip_0(\Omega)$.
From the estimates in~\eqref{e.equiv_moved} it follows that
inequality \eqref{e.pw_beta_hardy_natural} 
is comparable to \eqref{e.pointwise}, with $\kappa=2$
and $g=\lvert \nabla u\rvert$,
and therefore a pointwise $(p,w)$-Hardy inequality holds in $\Omega$ by~\eqref{e.minimality}.
Hence all assumptions of Theorem~\ref{t.abstract} are valid
and the claim follows from the $(p,w)$-Hardy inequality in Theorem~\ref{t.abstract},
applied with the bounded upper gradient $g_u$ 
that is given in connection with \eqref{e.minimality}.
\end{proof}

\begin{remark}
It is possible to extend Lemma~\ref{lem.w_is_poinc} and Theorem~\ref{t.euclidean}
also to some $-n<\beta<0$. 
In this case it is natural to add the condition that $w=0$ in $\Omega^c$.
The obstruction with $\beta<0$ is that clearly the
last inequality in \eqref{e.equiv_moved} is not valid for every $\beta<0$
if the ball $B(x,r)$ intersects the boundary of $\Omega$, since for small enough $\beta$
the integral in \eqref{e.equiv_moved} becomes infinite.
On the other hand, if the last inequality in \eqref{e.equiv_moved} is valid for some $\beta<0$, then
everything else in Lemma~\ref{lem.w_is_poinc} and Theorem~\ref{t.euclidean} works, 
and we conclude that for such $\beta<0$ the weighted pointwise $(p,\beta)$-Hardy inequality~\eqref{e.pw_beta_hardy_natural}
implies the weighted $(p,\beta)$-Hardy inequality~\eqref{e.beta_hardy}.

The validity of the last inequality in \eqref{e.equiv_moved} is closely related to the
Assouad dimension of $\partial\Omega$ via the so-called Aikawa condition, but
we omit any further discussion related to these concepts
and refer to~\cite{DydaEtAl2019,LehrbackTuominen2013} for details.
%
\end{remark}

\bibliographystyle{abbrv}
\def\cprime{$'$}

\end{document}